\documentclass[11pt]{amsart}

\usepackage[margin=2.5cm]{geometry}

\usepackage{epic}

\usepackage{enumitem}

\usepackage{amsmath,amsthm,amsfonts,amssymb,mathtools,mathrsfs,bm}

\usepackage{accents}

\usepackage{nicefrac}

\usepackage{tikz-cd}
\usetikzlibrary{babel}
\usetikzlibrary{decorations.pathmorphing}

\usepackage[utf8]{inputenc}

\usepackage[T1]{fontenc}

 \usepackage{lmodern}

\usepackage{float}

\usepackage{xcolor}
\definecolor{colordelink}{rgb}{0,0,0.50}
\definecolor{colordecite}{rgb}{0,0.5,0}
\definecolor{colordeurl}{rgb}{0,0.41,0.5}

\usepackage[pagebackref=true]{hyperref}
\hypersetup{
    colorlinks=true,
    linkcolor=colordelink,
    citecolor=colordecite,
	 urlcolor=colordeurl,
	 final
}

\usepackage[capitalize,nameinlink,noabbrev]{cleveref}

\usepackage{todonotes}

\usepackage{multirow}
\usepackage{cellspace}
\def\im{\operatorname{im}}
\def\dim{\operatorname{dim}}

\def\rk{\operatorname{rank}}

\def\id{\operatorname{id}}

\newcommand{\RR}{\mathbb{R}}

\newcommand{\CC}{\mathbb{C}}

\newcommand{\ZZ}{\mathbb{Z}}
\newcommand{\R}{\mathbb{R}}

\newcommand{\C}{\mathbb{C}}

\newcommand{\Z}{\mathbb{Z}}

\newcommand{\eqA}{\mathscr{A}}

\newcommand{\Ascr}{\mathscr{A}}

\newcommand{\Cbb}{\mathbb{C}}

\newcommand{\Sbb}{\mathbb{S}}
\newcommand{\Tbb}{\mathbb{T}}

\newcommand{\Ccal}{\mathcal{C}}

\newcommand{\Fcal}{\mathcal{F}}

\newcommand{\Lcal}{\mathcal{L}}

\newcommand{\Ocal}{\mathcal{O}}

\newcommand{\Rcal}{\mathcal{R}}
\newcommand{\Scal}{\mathcal{S}}

\newcommand{\Wcal}{\mathcal{W}}

\makeatletter
\newcommand{\tpitchfork}{%
  \vbox{
    \baselineskip\z@skip
    \lineskip-.52ex
    \lineskiplimit\maxdimen
    \m@th
    \ialign{##\crcr\hidewidth\smash{$-$}\hidewidth\crcr$\pitchfork$\crcr}
  }%
}
\makeatother

\theoremstyle{plain}
\newtheorem{theorem}{Theorem}[section]
\newtheorem{lemma}[theorem]{Lemma}
\newtheorem{corollary}[theorem]{Corollary}
\newtheorem{proposition}[theorem]{Proposition}

\theoremstyle{definition}
\newtheorem{definition}[theorem]{Definition}
\newtheorem{proposition/definition}[theorem]{Proposition/Definition}
\newtheorem{conjecture}[theorem]{Conjecture}
\newtheorem{example}[theorem]{Example}

\newtheorem{remark}[theorem]{Remark}
\newtheorem{notation}[theorem]{Notation}
\Crefname{notation}{Notation}{Notations}

\theoremstyle{remark}
\newtheorem*{note}{Note}

\usepackage{stackengine}
\def\rlwd{.4pt}
\def\rlht{1.1pt}
\def\shatvrule{\rule{\rlwd}{\rlht}}
\def\shat#1{%
  \setbox0=\hbox{$#1$}%
  \stackon[0pt]{\stackon[1pt]{\ensuremath{#1}}{%
    \shatvrule\kern\wd0\kern-\rlwd\kern-\rlwd\shatvrule}}%
    {\rule{\wd0}{\rlwd}}%
}

\newcommand{ \lk }{ \textnormal{lk} }

\def\C{\mathbb C}

\def\R{\mathbb R}
\def\Z{\mathbb Z}

\def\im{{\rm Im}}

\newcommand{\labelpar}{\label}

\author{R. Giménez Conejero and
Gerg\H{o} Pintér}

\title{Signature of the Milnor fiber of parametrized surfaces}

\address{
Alfr\'ed R\'enyi Institute of Mathematics, Re\'altanoda utca 13-15,
H-1053 Budapest, 
Hungary
}
\email{Roberto.Gimenez@uv.es}
\address{Department of Theoretical Physics, Institute of Physics,
Budapest University of Technology and Economics
\\
M\H{u}egyetem rkp. 3., H-1111 Budapest, Hungary
}
\email{pinter.gergo@ttk.bme.hu}

\subjclass[2020]{Primary 32S25, 32S50; Secondary 57R42} \keywords{Signature of 4-manifolds, Milnor fiber, deformations of map germs, complex surfaces, regular homotopy, Smale invariant, topological invariants}

\date{}

\begin{document}

\begin{abstract}
    We compute the signature of the Milnor fiber of certain type
of non-isolated complex surface singularities, namely, images of finitely determined holomorphic germs. An explicit formula  is given in algebraic terms. As a corollary we show that the signature of the Milnor fiber is a topological invariant for these singularities. The proof combines complex analytic and smooth topological techniques. The main tools are Thom--Mather theory of map germs and the Ekholm--Sz\H{u}cs--Takase--Saeki formula for immersions. We give a table with many examples for which the signature is computed using our formula.
\end{abstract}

\maketitle


\section{Introduction}

\subsection{Summary of the result}   In this paper, we combine several techniques to provide a formula for the signature of the Milnor fiber of certain type of non-isolated surface singularities. The singularities $(\Scal, 0) \subset (\C^3, 0)$ we consider are the images of holomorphic map germs $f: (\C^2, 0) \to (\C^3, 0)$ that have \textsl{mild non-isolated singularity} (i.e., finitely-determined germs). Recall that the Milnor fiber of a singularity is a nearby smooth fiber of its equation.

\begin{theorem}[see \cref{thm:main}]\label{th:main}
    Consider a finitely determined holomorphic germ $f: (\C^2, 0) \to (\C^3, 0)$. Let $(\Scal, 0)=\big(\im(f), 0\big)$ be its image and $F_f$ its Milnor fiber. Then, its signature $\sigma(F_f)$ satisfies
    \begin{equation*}
        \sigma (F_f)=\sigma(X)+T(f)-C(f).
    \end{equation*}
\end{theorem} 

In this theorem, $X$ is a manifold that we construct and $\sigma(X)$ is its signature, which can be computed by the explicit construction of $X$. The triple point number $T(f)$ and the cross-cap number $C(f)$  can be computed algebraically \cite{Mond1985, Mond1987}. In conclusion, this gives a computable way of expressing $\sigma(F_f)$. Indeed, we give a table with many examples using this formula in \cref{table:examples}.
\newline

Since all the ingredients of the formula are topological (by \cite{Bobadilla2019a} and \cite{Pinter2023}), we prove the following corollary.
\begin{corollary}[see \cref{cor:signature top}]
    The signature $\sigma(F_f)$ of \cref{th:main} is a topological invariant of map germs $f$.
\end{corollary}


The proof of \cref{th:main} is a mixture of complex analytic and smooth topological techniques. 

On the one hand, we deform $f$ into a more suitable map: we take a complex stable deformation in the sense of Thom-Mather theory (see, e.g., \cite{Mond2020}) and then we change to the smooth category to further deform this in the $\Ccal^\infty$ sense, by \cite{Nemethi2015}. On the other hand, we need a description of the boundary $\partial F_f$, given in \cite{Nemethi2018} (cf. \cite{Nemethi2012,Siersma1991a}). Then, using the works \cite{Nemethi2015,Nemethi2018}, we construct a smooth manifold $X$ together with a map (see these steps in \cref{fig:4figs} and $X$ in \cref{fig:X2d,fig:X3d}). The explicit construction of $X$ and the map is a big part of this work. 

Finally, we use the Ekholm--Sz\H{u}cs--Takase--Saeki formula \cite{Saeki2002}, which provides a regular homotopy invariant for immersions $N^3 \looparrowright \R^5$ in terms of the properties of a \emph{slice singular manifold} of the immersion.
We apply this formula for the inclusion $\partial F_f \hookrightarrow \Sbb^5 \subset \C^3$ using two different slice singular manifolds. The first one is the inclusion of the Milnor fiber, $F_f \hookrightarrow B^6 \subset \C^3$, the second one is the map from $X$ that we construct. We prove the theorem by comparing the invariants of both slice singular manifolds.
\newline

Applying the Ekholm--Sz\H{u}cs--Takase--Saeki formula to compute the signature of a Milnor fiber is an original idea of Andr\'{a}s N\'{e}methi, in discusions related to his joint work with the second author of this paper. This work is the first instance where this general plan is fulfilled, for this kind of singularities, and also one of the first applications of the formula for complex singularities.

\subsection{The context of the result} 
The Milnor fiber plays a central role in the study of local singularities. Considering a hypersurface singularity $(\mathcal{S}, 0) =\big(g^{-1}(0), 0\big) \subset (\C^{n+1}, 0)$, recall that its Milnor fiber, defined as $F\coloneqq g^{-1}(\delta) \cap B^{2n+2}_{\epsilon}$ for $0<|\delta|\ll|\epsilon|\ll 1$, is a smooth oriented manifold of real dimension $2n$.

In the case of isolated hypersurface singularities (i.e., if $dg^{-1}(0) \cap \Scal \subset \{0\}$) the Milnor fiber is rather well understood. It has the homotopy type of a wedge
of $n$-spheres, whose number is called the Milnor number $\mu$, that is, $\mu=|\chi(F)-1|$ (where $\chi(F)$ is the Euler characteristic of the Milnor fiber). Moreover, $\mu$ is also the codimension of the Jacobian ideal, generated by the partial derivatives of $g$. Also, for isolated hypersurface singularities the boundary of the Milnor fiber
is diffeomorphic to the link $K \coloneqq \mathcal{S} \cap \Sbb_{\epsilon}^{2n+1}$ of the singular germ, and also with the boundary
of any resolution of the germ \cite{Milnor1968,Nemethi1999}. This coincidence produces several nice formulas connecting the invariants of
these fillings, formulas of Laufer \cite{Laufer1977} or Durfee \cite{Durfee1978} and their generalizations, see e.g. \cite{Wahl1981}, and the Durfee conjecture \cite{Kollar2017}. 

For non--isolated hypersurface singularities in $(\C^3,0)$ the situation is
more complicated. The Milnor fiber is not necessarily simply connected and the Jacobian ideal has infinite codimension. The link $K$ of $(\Scal, 0)$ is not smooth, hence it cannot be diffeomorphic to the boundary of the Milnor fiber.
A general algorithm to construct the boundary of the Milnor fiber of any non-isolated surface singularity $ (\Scal, 0) \subset ( \C^3, 0) $ is presented in \cite{Nemethi2012}, although it is rather technical. For particular families
of singularities there are more direct descriptions of $\partial F$ from the peculiar intrinsic geometry of the germ, see e.g. \cite{Michel2007} or \cite{Bobadilla2014}, and also \cite{Nemethi2018} in the context we study.

For a finitely determined map germ $f: (\C^2, 0) \to (\C^3, 0)$ there are two different candidates for the generalization of the Milnor number. An analogue of the Milnor fiber is the image of a stable deformation of $f$ (called \textsl{disentanglement}), but it is not a smooth manifold. Indeed,  it contains the cross-caps, the triple values (whose number is $C(f)$ and $T(f)$, respectively), and the curve of double values. Even so, the disentanglement is homotopy equivalent to a wedge of 2-spheres, whose number is called the \emph{image Milnor number} $\mu_I(f)$, see \cite{Mond1991}. It can be expressed as (see \cite{Marar1989,Mond1991})
\begin{equation}\label{eq:imagemilnor}
    \mu_I(f)=\frac{1}{2}\big(\mu (D) +C(f)-4 T(f)-1\big),
\end{equation}
where $\mu(D)$ is the Milnor number of the double point set $(D,0) \subset (\C^2, 0)$.

The other candidate for the Milnor number is the second Betti number of the Milnor fiber $b_2(F_f)$.
 In our context the first Betti number is zero, hence $b_2(F_f)=\chi(F_f)-1$ and it is expressed as (see \cite{Bobadilla2019a}, cf. \cite{Siersma1991,Massey1992})
\begin{equation}\label{eq:bettitwo}
    b_2(F_f)=\mu(D)+2C(f)-3T(f)-1.
\end{equation}

The contribution of this paper is \cref{th:main}, which provides a similar formula for the signature of the Milnor fiber in our context. In general, the signature of an oriented 4-manifold is an invariant by homeomorphisms. Conversely,  together with the parity and rank, it determines the intersection form of 4-manifolds that are closed and simply-connected, which \textsl{almost} characterizes them modulo homotopy and homeomorphism by the Milnor-Whitehead theorem and Freedman's theorem (see \cite[Theorems 1.2.25 and 1.2.27]{Gompf1999}, \cite{Freedman1982,Freedman1990}). The signature is also an important invariant of the Milnor fiber of isolated complex singularities. It induces a finer characterisation of $F_f$ together with the Milnor number $\mu$ (which determines the homotopy type of $F_f$).
The signature appears in the Durfee type formulas mentioned before.
\newline

In general, it is unknown whether two ambient-homeomorphic hypersurface singularities have homeomorphic Milnor fiber or not, except for isolated singularities where this holds (see \cite{Saeki1989}, cf. \cite{King1978,Perron1985}). By \cite{Le1973}, the \emph{homotopy type} of the Milnor fibre is a topological invariant. In our case of non-isolated singularities, i.e. for images  of finitely determined germs, the boundary $\partial F_f$ is topological \cite{Pinter2023}. Now, for these singularities, we show that the signature $\sigma (F_f)$ is also a topological invariant. Hence, it is natural to conjecture the following.

\begin{conjecture}
Two ambient-homeomorphic hypersurface singularities have homeomorphic Milnor fibers.
\end{conjecture}

Note that, while the previously known formulas like \cref{eq:bettitwo,eq:imagemilnor} are proved in the complex analytic category, our proof for \cref{th:main} escapes to the real $\mathcal{C}^{\infty}$ world, namely, immersion theory. Immersion-theoretical approaches for problems of complex singularities appear in several works, e.g., in \cite{Nemethi2015, Pinter2023a, Pinter2023} the associated immersion of finitely determined germs plays the key role, and in \cite{Katanaga2014,Ekholm2006}
other type of singularities are investigated. Interestingly, \cite{Ekholm2006} applies an Ekholm--Sz\H{u}cs--Takase--Saeki type formula in reverse compared to us. Namely, it identifies the inclusions of the Brieskorn exotic 7-spheres in $\Sbb^9$ up to regular homotopy in terms of the signature of the Milnor fiber, while we identify the latter invariant from a regular homotopy invariant of the boundary.

\subsection{Structure of the article}

In \cref{sec:prel} we introduce all the ingredients that are necessary to understand this work. We start with properties of finitely determined holomorphic germs. Then, we summarize the Ekholm--Sz\H{u}cs--Takase--Saeki formula from immersion theory, and its application for the associated immersion of complex germs.

In \cref{s:ssmftbotmf} we construct a slice singular manifold for the boundary of the Milnor fiber, that is necessary for the proof of our main theorem. 

\cref{sec:main} gives the proof of main theorem, which is based on the application of the Ekholm--Sz\H{u}cs--Takase--Saeki formula for the slice singular manifold we constructed.

In \cref{s:ex} we compute the signature of the Milnor fiber for several families of examples.  This is summarized in \cref{table:examples}.

\subsection{Acknowledgments}

We thank Andr\'as N\'emethi for suggesting the topic of this paper and numerous fruitful discussions. We are grateful to Andr\'as Sz\H{u}cs,  Marco Marengon,  Andr\'as S\'andor for answering our questions and sharing their knowledge on related topics. GP thanks his physicist colleagues Andr\'as P\'alyi, Gy\"orgy Frank, Zolt\'{a}n Guba, D\'aniel Varjas and J\'anos Asb\'oth for the new inspiration to the singularity theory research.

GP was supported by the Ministry of Culture and Innovation and the National Research, Development and Innovation Office within the Quantum Information National Laboratory of Hungary (Grant No. 2022-2.1.1-NL-2022-00004) as well as by the National Research, Development and Innovation Office via the OTKA Grant No. 132146.

\section{Preliminaries}\label{sec:prel}

\subsection{Map germs and invariants}\label{sec: intro germs}

Here we give a quick introduction to Thom-Mather theory of holomorphic map germs, the reader is referred to \cite[Section 1.2]{Robertothesis} for an in-deep introduction and \cite{Mond2020} for a general modern reference on this topic.

We are interested in map germs $f:(\CC^2,0)\to(\CC^3,0)$ modulo $\Ascr$\textsl{-equivalence}. More precisely, two map germs $f,f'$ are $\Ascr$\textit{-equivalent} (or \textit{left-right equivalent}) if there are biholomorphisms $\varphi$ and $\psi$ that make the following diagram commutative
\begin{equation}\label{eq:defAeq}
\begin{tikzcd}
 (\CC^2,0) \arrow[r, "f" ]\arrow[d, " \varphi"',"\sim"  {anchor= north, rotate=90, inner sep=.6mm}]& (\CC^3,0)\arrow[d, " \psi","\sim"'  {anchor= south, rotate=90, inner sep=.3mm}] \\
 (\CC^2,0) \arrow[r, "f'" ]& (\CC^3,0)
\end{tikzcd}.
\end{equation}
We say that a map germ $f$ is \textit{stable} if $f$ is $\Ascr$-equivalent to any of its deformations $f_t$ (i.e., the germ given by $f_t$ at some point), and it is \textit{unstable} otherwise.

\begin{example}[Complex Whitney umbrella, cross-cap]\label{ex:cc}
The map germ $f:(\CC^2,0)\to(\CC^3,0)$ given by $f(u,v)=(u,v^2,uv)$ is stable. In particular, the following deformation $f_t$ is $\Ascr$-equivalent to $f$ at the origin:
\[\begin{tikzcd}[ampersand replacement = \&]
	\&[-25pt] \textcolor[rgb]{0.32,0.32,0.32}{(u,v)} \&[-40pt] \&[-40pt]  \&[-60pt] \textcolor[rgb]{0.32,0.32,0.32}{(u,v^2,uv)} \&[-45pt]  \\[-15pt]
	\textcolor[rgb]{0.32,0.32,0.32}{(u+tu^2,v)} \& {(\mathbb{C}^2,0)} \& \&  \& {(\mathbb{C}^3,0)} \&[-30pt] \textcolor[rgb]{0.32,0.32,0.32}{(x+tx^2,y,z)} \\[-20pt]
      \& \& \& \textcolor[rgb]{0.6,0.6,0.6}{\big(u+tu^2,v^2,(u+tu^2)v\big)}\&  \&\\
     \& \& \textcolor[rgb]{0.6,0.6,0.6}{(u,v)} \& \& \&    \\[-20pt]
	\textcolor[rgb]{0.32,0.32,0.32}{(u,v)} \& {(\mathbb{C}^2,0)} \& \&  \& {(\mathbb{C}^3,0)} \& \textcolor[rgb]{0.32,0.32,0.32}{(x,y,z)} \\[-15pt]
	\& \textcolor[rgb]{0.32,0.32,0.32}{(u,v)} \& \& \& \textcolor[rgb]{0.32,0.32,0.32}{(u,v^2,uv+tu^2v)} \& 
	\arrow["f", from=2-2, to=2-5]
	\arrow[maps to, color={rgb,1:red,0.32;green,0.32;blue,0.32},from=1-2, to=1-5]
	\arrow["\varphi","\sim" {anchor= north, rotate=90, inner sep=.6mm}, from=5-2, to=2-2]
	\arrow["\psi"',"\sim" {anchor= south, rotate=90, inner sep=.6mm}, from=5-5, to=2-5]
	\arrow["{f_t}"', from=5-2, to=5-5]
	\arrow[maps to, color={rgb,1:red,0.32;green,0.32;blue,0.32},from=6-2, to=6-5]
	\arrow[maps to, color={rgb,1:red,0.32;green,0.32;blue,0.32},from=5-6, to=2-6]
	\arrow[maps to, color={rgb,1:red,0.32;green,0.32;blue,0.32},from=5-1, to=2-1]
   \arrow[maps to,color={rgb,1:red,0.6;green,0.6;blue,0.6}, dashed, from=4-3, to=3-4]
\end{tikzcd}.\]
\end{example}

\begin{figure}[htbp]
	\centering
		\includegraphics[width=0.4\textwidth]{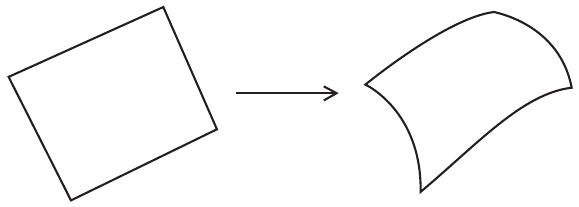}\quad\quad\includegraphics[width=0.4\textwidth]{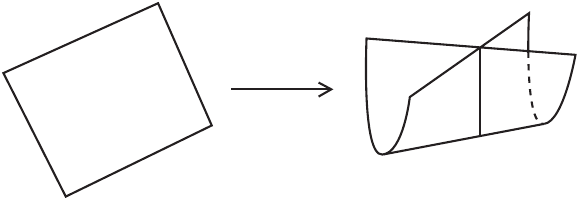}\\
		\includegraphics[width=0.4\textwidth]{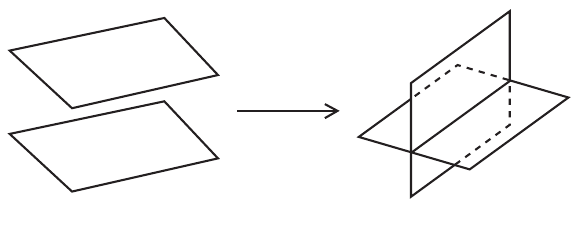}\quad\quad\includegraphics[width=0.4\textwidth]{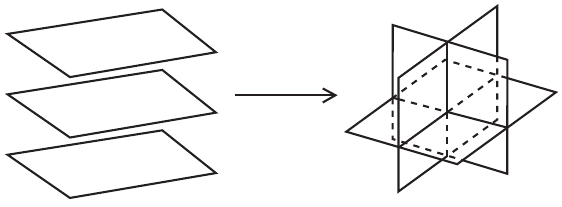}      
	\caption{Stable singularities from $\CC^2$ to $\CC^3$.}
	\label{fig:stables}
\end{figure}

Indeed, the only stable singularities (multigerms) from $\Cbb^2$ to $\Cbb^3$ are \textit{injective immersions}, \textit{Whitney umbrellas} or \textit{cross-caps} (of dimension zero, as the germ given in \cref{ex:cc}), \textit{transverse double points} (transverse intersection of two regular branches, of dimension one), and \textit{transverse triple points} (normal crossing intersection of three regular branches, of dimension zero); see \cref{fig:stables}.

There is a kind of map germs where the question of whether or not two maps are $\Ascr$-equivalent is simpler than in the general case.

\begin{definition}\hypertarget{htaf}{}
A map germ $f:(\CC^2,0)\to(\CC^3,0)$ is \textit{$k$-determined} if it is $\Ascr$-equivalent to $f'$ whenever their Taylor polynomial of order $k$ at $0$ coincide. If a map germ is $k$-determined for some $k$ we say that it is \textit{finitely determined} or \textit{$\Ascr$-finite}.
\end{definition}

By Mather-Gaffney criterion (see \cite[Theorem 4.5]{Mond2020}), $\Ascr$-finite map germs are characterized by having isolated instability, i.e., every germ away from the origin is stable. For map germs $f:(\CC^2,0)\to(\CC^3,0)$, this is equivalent  to $f$ being a stable (i.e., transverse) immersion away from the origin for a small enough representative, to avoid triple points and Whitney umbrellas. 

The image of a finitely determined germ $f:(\C^2, 0) \to (\C^3, 0)$ is a  surface $(\mathcal{S}, 0)=\big(\im(f), 0\big)$ in $(\C^3, 0)$ with a non-isolated singularity (except when $f$ is regular, cf. \cite[
Theorem 1.3.1.]{Nemethi2015}). Its transverse type is $A_1$ (cf. \cref{lemma:TC} below), and its normalization map is $f$ itself.


Despite having non-isolated singularities, there is a way of taming these images: every $\eqA$-finite germ has a deformation $f_s$ that is stable at every point. So, in particular, $f_s$ has at most cross-caps, transverse double points and transverse triple points (see more details in \cite[Section 5.5]{Mond2020}). 

In this context, if we have two $\eqA$-finite germs $f$ and $f'$ that are $\eqA$-equivalent, their stable deformations are left-right equivalent as maps, therefore, the number of cross-caps $C(f)$ (sometimes $W(f)$ in the literature) and triple points $T(f)$ that appear are $\Ascr$-invariants of $f$,  introduced by Mond in \cite{Mond1985, Mond1987}. They appear in several different contexts, see for example \cite{Mond1991, Marar1989, Marar2012, Marar2014, Gergothesis}. Furthermore, there is an algebraic way of computing them.

\begin{lemma}[see {\cite[Exercise E.4.2 and Corollary 11.12]{Mond2020}}]\label{lemma:TC}
For an $\Ascr$-finite map germ $f:(\CC^2,0)\to(\CC^3,0)$,
$$\begin{aligned} 
C(f)&=\dim_\CC \frac{\Ocal_{\CC^2,0}}{\Rcal_f},\textnormal{ and}\\
 \quad T(f)&=\dim_\CC \frac{\Ocal_{\CC^3,0}}{\Fcal_2},
\end{aligned}$$
where $\mathcal{R}_f$ is the ramification ideal, generated by the $2\times 2$ minors of the Jacobian matrix $df$, and $\mathcal{F}_2$ is the second Fitting ideal.
\end{lemma}

\subsection{The associated immersion}\label{ss:associmm} 

\begin{figure}
	\centering
		\includegraphics[width=0.9\textwidth]{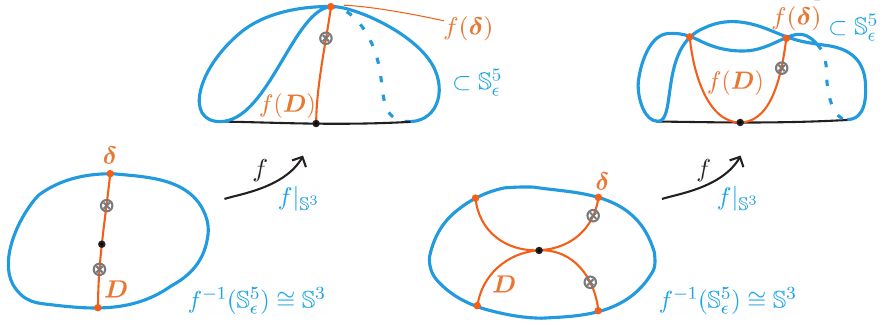}
	\caption{An example of a map with a twisted component (left, cf. \cref{ex:cc,ex:s0}) and another with an untwisted component (right, cf. $S_1$ in \cref{ex:sk}). The associated immersion \textcolor[rgb]{0,0.6,0.85}{$f|_{\Sbb^3}: f^{-1}(\Sbb^5_{\epsilon}) \cong \Sbb^3 \to \Sbb^5_{\epsilon}$}; the double point sets \textcolor[rgb]{0.95,0.44,0.13}{$\bm{D}$} and \textcolor[rgb]{0.95,0.44,0.13}{$f(\bm{D})$}; and the link \textcolor[rgb]{0.95,0.44,0.13}{$\bm{\delta}$} of \textcolor[rgb]{0.95,0.44,0.13}{$\bm{D}$}, coinciding with the double point set of \textcolor[rgb]{0,0.6,0.85}{$f|_{\Sbb^3}$} in the source.} 
	\label{fig:simpleEx}
\end{figure}

Every $\Ascr$-finite map germ $f: (\C^2, 0) \to (\C^3,0)$ induces a stable immersion $\Sbb^3\looparrowright \Sbb^5$. Indeed, the preimage $f^{-1}\big(\Sbb^5_\epsilon\big)\subset\CC^2$ of a small enough sphere $\Sbb^5_\epsilon\subseteq \Cbb^3$ centered at $0$, is diffeomorphic to $\Sbb^3$. The restriction of $f$ to $f^{-1}\big(\Sbb^5_\epsilon\big)$ is a stable immersion that we denote by $f|_{\Sbb^3}:\Sbb^3\looparrowright \Sbb^5$. Moreover, different choices in the definition of $f|_{\Sbb^3}$ give stable immersions $\Sbb^3\looparrowright \Sbb^5$ that are (regular) homotopic through \textit{stable} immersions. Note that stable in this context means that the immersion has at most transverse double points. See the details in \cite[Section 2.1]{Nemethi2015} or \cite[Subsection 1.1.2.]{Gergothesis}, and also \cref{fig:simpleEx}.

We can say even more: Since $\Ascr$-finite map germs are stable away from the origin, any (complex) deformation $f_t$ induces, by an identification of the source $f_t^{-1}(\Sbb^5_\epsilon)$, the same stable immersion $\Sbb^3\looparrowright \Sbb^5$ modulo (regular) homotopy through stable immersions. In particular, this holds for any stable deformation $f_s$. Cf. \cite[Proposition 9.1.2]{Nemethi2015}.

The image of the immersion $f|_{\Sbb^3}$ is the link of the image of $f$, that is,  $K_f \coloneqq \im (f) \cap \Sbb^5_{\epsilon} $.

\subsection{The image and the double points}\label{ss:imdoub} 

An interesting and important object is the \textsl{double point set} of $f$, denoted usually by $D(f)$, which hides a lot of information about the map germ $f$ (see, for example, \cite{GimenezConejero2022b} or \cite{Marar2012}). The \textit{double point set} in this context can be defined as the closure of the double point locus with reduced structure, i.e., the germ given by
$$ D(f)\coloneqq \textnormal{cl}\big\{x\in\CC^2: \exists y\neq x \text{ s.t. } f(x)=f(y)\big\}.$$
For simplicity, we will use $\bm{D}$ instead of $D(f)$.

Let $ g: ( \C^3, 0) \to ( \C, 0) $ be the reduced equation of $\big( \im(f),0\big)$. As we pointed out before, this is a germ of a complex surface with non-isolated singularity, whose singular locus coincides with $\big(f(\bm{D}),0\big)$. In particular, observe that $\bm{D}$ is a plane curve and $f(\bm{D})$ is a curve in $\CC^3$ considered with their reduced structure (see \cref{fig:simpleEx}).

Let $f(\bm{\delta})$ denote the link of $f(\bm{D})$, i.e., $f(\bm{\delta})\coloneqq f(\bm{D})\cap\Sbb^5_\epsilon$, and, by analogy, let $\bm{\delta}\coloneqq f^{-1}\big(f(\bm{\delta})\big)$, which is also the link of $\bm{D}$ since the preimage $f^{-1}(B_{\epsilon}^6) \cong B^4$ serves as a Milnor ball for $\bm{D}$; both considered with their natural orientation.
Furthermore, observe that $\bm{\delta}$ and $f(\bm{\delta})$ induce the double points and the double values of $f|_{\Sbb^3}$, as defined in \cref{ss:associmm}:
\[\begin{tikzcd}
	{(\Cbb^2,0)} & {(\Cbb^3,0)} \\
	{\Sbb^3} & {\Sbb^5} \\
	{\bm{\delta}} & {f(\bm{\delta})}
	\arrow["f", from=1-1, to=1-2]
	\arrow["{f|_{\Sbb^3}}", from=2-1, to=2-2]
	\arrow[dashed, hook, from=2-1, to=1-1]
	\arrow[dashed, hook, from=2-2, to=1-2]
	\arrow["{f|_{\Sbb^3}}", from=3-1, to=3-2]
	\arrow[hook, from=3-1, to=2-1]
	\arrow[hook, from=3-2, to=2-2]
\end{tikzcd},\]
where the dashed arrows only make sense after taking a representative of $f$.

Let $\bigsqcup_{i=1}^n \delta_i=\bm{\delta}$ be the decomposition of $\bm{\delta}$ corresponding to the irreducible components $D_i$ of $\bm{D}$. Note that $\delta_i\cong\Sbb^1$ for every $i$. Some of these components $\delta_i$ are mapped one-to-one to their images $f(\delta_i)$, we call those \textit{untwisted components}. As they are part of the double points, the untwisted components are mapped in pairs to the same connected component $f(\delta_i)$ of $f(\bm{\delta})$. Those components $\delta_i$ that are mapped two-to-one to their image $f(\delta_i)$ are called \textit{twisted components}. Hence, as $f(\delta_i)\cong\Sbb^1$ in both cases, untwisted components provide a trivial double cover of $\Sbb^1$ and twisted components a non-trivial double cover (see \cref{fig:simpleEx}). 
\begin{notation}\label{not:fi}
For simplicity, we use the index notation for the components of $f(\bm{D})$
$$ f\{i\}\coloneqq\begin{cases}
\left\{i,i'\right\} & \text{if }\delta_i\neq\delta_{i'}\text{ are untwisted and }f(\delta_i)=f(\delta_{i'}),\\
\left\{i\right\} & \text{if }\delta_i\text{ is twisted.}
\end{cases} $$
\end{notation}

We also need to know the topology around the different $f(\delta_i)$ depending on whether $\delta_i$ is twisted or untwisted.
\begin{lemma}[cf. {\cite[Section 4.2]{Nemethi2018}}]\label{lem:twORuntw}
Let $T_{f\{i\}} \cong \Sbb^1 \times B^4 $ be a  small enough closed tubular neighborhood of $f(\delta_i)$ in $\Sbb_\epsilon^5$.
Then $T_{f\{ i\}}\cap K_f$ is the total space of a locally trivial fiber bundle given by the retraction $r$ to $f(\delta_i)\cong\Sbb^1$,
\[\begin{tikzcd}
	{B^2\vee B^2} & T_{f\{ i\}}\cap K_f \\
	& {f(\delta_i)}
	\arrow["r", from=1-2, to=2-2]
	\arrow[hook, from=1-1, to=1-2]
\end{tikzcd}\]
such that a geometric monodromy is trivial if $\delta_i$ is untwisted and it swaps the discs in $B^2\vee B^2$ if it is twisted.
\end{lemma}
\begin{proof}
Recall that the  intersection of $\big(\im(f),0\big)$ with a complex plane  transverse to $f(\bm{D})$ at points in $f(\bm{D})\setminus\left\{0\right\}$ is an $A_1$ singularity $\left\{x'y'=0\right\}\subset \CC^2$.
Hence, around a point in $f(\bm{D}) \setminus \{0\}$, the image of $f$ is isomorphic to the variety of equation $xy=0$ in $\CC^3$, showing that the fiber is $B^2\vee B^2$.
Then, its intersection with $\Sbb^5$ is homeomorphic to $(B^2\vee B^2) \times I$ for some interval $I$. This proves that $r$ defines a locally trivial fiber bundle.

The isotopy class of the geometric monodromy is given by the mapping $f$, which gives $T_{f\{ i\}}\cap K_f$ as the image of two disconnected copies of $ B^2\times\Sbb^1$ in the untwisted case and as the image of only one copy of $ B^2\times\Sbb^1$ in the twisted case.
\end{proof}

Notice that the analytic version of \cref{lem:twORuntw} (which deals only with the topology) is proved in \cite[Section 4.2]{Nemethi2018}: There is (locally) a splitting $g=g_1\cdot g_2$ of $\im(f)$ along $f(D_i)\setminus\left\{0\right\}$  that glues globally in an analytic way, permuting of components in the twisted case.

\subsection{The boundary of the Milnor fiber} \label{ss:surgerypartialf}
\subsubsection{The surgery providing $\partial F_f$}\label{sss:surgerypartialf part1}
Recall that the Milnor fiber $F_f$ of $\big(\im(f),0\big)$ in $(\CC^3,0)$ is given by $F_f= g^{-1} ( \delta ) \cap B^{6}_{ \epsilon }$ for some $0<|\delta|\ll\epsilon\ll1$, where $g$ is the reduced equation of $\im(f)$. The Milnor fiber $F_f$ is a smooth oriented 4-manifold with boundary $\partial F_f \subset \Sbb^5_{ \epsilon }$. Since $\big(\im(f),0\big)$ is a surface with a non-isolated singularity, its link $K_f = \im(f) \cap \Sbb^5_{ \epsilon }$ is not smooth and, hence, there is obviously a change from $K_f$ to $\partial F_f$. 

Furthermore, $K_f$ is the image of $f|_{\Sbb^3}:\Sbb^3\looparrowright \Sbb^5$ given in \cref{ss:associmm}. Now we summarize the algorithm from \cite{Nemethi2018} to construct $\partial F_f$ from $f^{-1}(K_f)\cong \Sbb^3$ as a (generalized) surgery  along $\bm{\delta}$. See also \cite[Chapter 5]{Gergothesis}. The construction given in \cref{th:milnperpart1} below is standard, see, for example, \cite{Nemethi2012} or \cite{Siersma1991a}, however the hard part is the computation of the gluing coefficients called \emph{vertical indices}, summarized in \cref{sss:vertind}.

The construction is based on the following observation. $\partial F_f$ decomposes as the union of $\partial F_f \setminus \mathring{\bm{T}}$ and $\partial F_f \cap \bm{T}$, where $\bm{T}=\bigsqcup_{f\{i\} } T_{f\{i\}}$ is a tubular neighborhood of $f(\bm{\delta})$, see \cref{lem:twORuntw}, and $\mathring{\bm{T}}$ denotes its interior. 

The first part, $\partial F_f \setminus \mathring{\bm{T}}$ is diffeomorphic to $K_f \setminus \mathring{\bm{T}}$, moreover, they can be identified by an isotopy in $\Sbb^5 \setminus \mathring{\bm{T}}$. 

The second part is a locally trivial fiber bundle over $f(\delta)$ (cf. \cref{lem:twORuntw}),
\begin{equation}\label{eq:TC2}
    \begin{tikzcd}
	{\Sbb^1 \times [-1,1]} & T_{f\{ i\}}\cap F_f \\
	& {f(\delta_i)}
	\arrow[, from=1-2, to=2-2]
	\arrow[hook, from=1-1, to=1-2]
\end{tikzcd}
\end{equation}
whose fiber is the Milnor fiber of the transverse $A_1$ singularity. More precisely, we take a complex plane in $\CC^3$ that is transverse to $f(\bm{\delta})$ at a point, this gives an $A_1$ singularity in $\CC^2$ whose Milnor fiber is the cylinder $\{x'y'=\delta\} $.
The geometric monodromy of this bundle depends on whether the component $f(\delta_i)$ is untwisted or twisted (see \cite[Proof of Proposition 4
1]{Nemethi2018}).  

For untwisted components the geometric monodromy is trivial, implying that $T_{f\{ i\}}\cap F_f $ is diffeomorphic to $\Sbb^1 \times \Sbb^1 \times [-1,1]$ if $f\{i\}=\{i,i'\}$ is untwisted. However, for a twisted component $f\{i\}=\{i\}$, the geometric monodromy is given by rotating a cylinder $\Sbb^1\times\left[-1,1\right]$ half a twist, $(z,t)\mapsto(\overline{z},-t)$. Here $\Sbb^1 $ is considered as the unit sphere in $\C$, and $\overline{z}$ denotes the complex conjugate. Hence, in the twisted case, $T_{f\{ i\}}\cap F_f $ is diffeomorphic to the
3-manifold with boundary
\begin{equation}
 Z_\partial\coloneqq  \frac{\left[-\pi,\pi\right]\times\Sbb^1\times\left[-1,1\right]}{(-\pi,z,t)\sim(\pi,\overline{z},-t)}.
\label{eq:Y}
\end{equation}
It obviously fits into the following locally trivial fibration that has geometric monodromy $(z,t)\mapsto(\overline{z},-t)$,
\begin{equation}\label{eq:Zpartfiber}
\begin{tikzcd}
	{\Sbb^1\times\left[-1,1\right]} & Z_\partial\\
	& {\Sbb^1}
	\arrow["P_1", from=1-2, to=2-2]
	\arrow[hook, from=1-1, to=1-2]
\end{tikzcd},\end{equation}
where $P_1$ is the projection on the first component.  Furthermore, it is not hard to see that the boundary $\partial Z_\partial$ is a torus,  see \cite[Section 3]{Nemethi2018}, where more details about the geometry of $Z_{\partial}$ is given. The choice of 
the notation $Z_\partial$ will be clear in \cref{ss:4mfld}.

The pieces $\partial F_f \setminus \mathring{\bm{T}}$ and $\partial F_f \cap \bm{T}$ glue together along the torus boundaries corresponding to the components of $f(\bm{\delta})$. This gluing can be realized in the source of $f$, instead of the target, as follows. Define $\bm{N} \coloneqq f^{-1}(\bm{T}) \subset \Sbb^3$. It decomposes as $\bm{N}=\bigsqcup_i N_i$, where $N_i$ is  a closed tubular neighborhood of $\delta_i \subset \Sbb^3$, diffeomorphic to $B^2 \times \Sbb^1$. Obviously, $f^{-1}(\Sbb^5 \setminus \mathring{\bm{T}})=\Sbb^3 \setminus \mathring{\bm{N}}$. We define the gluing pieces
\begin{equation*}
    Y_{f\{i\}} :=
    \left\{
\begin{array}{cc}
\Sbb^1 \times \Sbb^1 \times [-1,1] &  \textnormal{ if } f\{i\}=\{i, i'\} \\

Z_{\partial} & \textnormal{ if } f\{i\}=\{i\}.
\end{array}
    \right.
\end{equation*}

Based on the above argument, the Milnor fiber is given by the following result, see \cite[Proposition 4.1]{Nemethi2018}.
\begin{theorem}\label{th:milnperpart1}  One has an orientation-preserving diffeomorphism
\begin{equation}\label{eq:rag}  
   \partial F \cong \left[ \big(\Sbb^3 \setminus\mathring{\bm{N}}  \big)\sqcup \bigsqcup_{f\{i\}} Y_{f\{i\}}   \right] 
	\Bigg / \Phi_\partial \textnormal{,} 
   \end{equation}
	where  the gluing 
 \begin{equation*} \Phi_\partial: \partial \left(\Sbb^3 \setminus\bigsqcup_{i} \mathring{N}_i  \right) \to -\partial \left( \bigsqcup_{f\{i\}} Y_{f\{i\}} \right)
 \end{equation*}
 is induced by a
	collection $\Phi^{f\{i\}}_\partial$ of diffeomorphisms
	\[ \Phi_\partial^{f\{i\}}:\begin{cases}
\partial N_i \sqcup \partial N_{i'} \to \partial (\Sbb^1 \times \Sbb^1 \times [-1,1]) & \text{if }f\{i\}=\left\{i,i'\right\}\\
 \partial N_i \to  \partial Z_{\partial} & \text{if }f\{i\}=\left\{i\right\}.
\end{cases} \]
\end{theorem}

We have to determine the gluing maps $\Phi_\partial^{f\{i\}}$. In order to do so, first observe that in the untwisted case the gluing map can be simplified to an identification $\partial N_i \to -\partial N_{i'}$. Hence in both cases we have a diffeomorphism between tori. In order to describe it we need to fix a homological trivialization.

\begin{definition}
Given a solid torus $\Sbb^1\times B^2\subset\Sbb^3$ that retracts to the knot $\delta$, and its boundary $\Tbb^2$, an \textit{oriented meridian} of $\Tbb^2$ is a closed curve whose linking number in $\Sbb^3$ with $\delta$ is $1$ and bounds a disc in $\Sbb^1\times B^2$, and a \textit{topological longitude} (or \textit{Seifert framing of} $\delta$) is a closed curve that generates $ H_1 (\Sbb^1\times B^2, \Z )$, whose linking number with $\delta$ is $0$ and has the same orientation as $\delta$.
\end{definition} 

We choose the generators of $H_1( \partial N_i, \Z)$ induced by the oriented meridian and the topological longitude. This is enough to describe the surgery on the untwisted parts. However, the torus $\partial Z_\partial\cong\Tbb^2$ is not embedded in $\Sbb^3$, so we need to take generators of its homology to proceed with the surgery. We take as \textsl{oriented longitude} any closed orbit  of any point in $\partial Z_\partial$ along the base $\Sbb^1$ in \cref{eq:Zpartfiber} (which takes two loops), and as \textsl{oriented meridian} any boundary component of any cylinder fiber. The homology classes of these cycles are independent of choices. See more details in \cite[Section 3]{Nemethi2018}. By construction, we have the following.

\begin{theorem}\label{th:milnperpar2}
    Using the (homological)  trivialization of the tori $ \partial N_i $ and $ \partial Z_{\partial} $ determined by the pair (meridian, longitude), each $\Phi_\partial^{f\{i\}}$ in \cref{th:milnperpart1} is given by a matrix    
	\begin{equation}\label{eq:milnpermatrix}
	\left( \begin{array}{cc}
	-1 & \mathfrak{vi}_{f\{i\}} \\
	0 & 1 \\
	\end{array} \right).
	\end{equation} 
\end{theorem}
   
Therefore, the gluing is determined by one integer for each component. 
\begin{definition}\label{def:verticalindex}
    The gluing coefficients $\mathfrak{vi}_{f\{i\}}$ appearing in \cref{th:milnperpar2} are called vertical indices (in \cite{Nemethi2018}, the $\mathfrak{vi}_{f\{i\}}$ are denoted by $\mathfrak{vi}_j$).
\end{definition}

\begin{remark}
Using \cref{th:milnperpart1,th:milnperpar2}, it is possible to present $ \partial F$ as a plumbed $3$-manifold, see \cite[Sections 4.6 and 4.7]{Nemethi2018}. It is very important to notice, however, that our construction of $X$ in \cref{ss:4mfld} below is not the 4-manifold corresponding to that plumbing graph. It seems that they are related by a sequence of blow-ups.
\end{remark}

\subsubsection{Computation of the vertical indices}\label{sss:vertind}
For the sake of completion, we present here a sketch of the computation of the vertical indices, but it is only relevant for calculations: they are implicit in our signature formula \cref{th:main} (see \cref{ss:signature}) and we use them explicitly in the computations of examples in \cref{s:ex}. The reader can ignore this subsection if the details are not needed.
\newline

 The computation of the vertical indices is special for these singularities and it is described now. One can also see this (formulated as a definition) in \cite[Definition 4.10]{Nemethi2018} via non-trivial constructions and statements, using an \textsl{aid germ} $\tau$ and Taylor expansions of $\tau$ and $g$.

\begin{definition}
A germ $ \tau: ( \C^3, 0) \to ( \C ,0) $, or $H\coloneqq \tau^{-1}(0)$, is called a \textit{transverse section along} $f(\bm{D}) $ if
\begin{enumerate}[label=(\roman*),font=\itshape]
	\item $ \big(f(\bm{D}), 0\big) \subset H$,
   \item $H$ is smooth at any point $p\in f(\bm{D}) \setminus \{ 0 \}$, and
   \item $H$ is transverse to both components of $\im(f)$ at any $p\in f(\bm{D}) \setminus \{ 0 \}$.
\end{enumerate}
\end{definition}

 Transverse section always exists, see \cite[Proposition 4.4]{Nemethi2018}. We fix one. Clearly, $ (\tau \circ f)^{-1}(0) \subset ( \C^2, 0) $ decomposes as $ \bm{D} \cup \bm{D}_{\sharp} $ for some curve $ \bm{D}_{\sharp} $ (not necessarily reduced).  Then, we define
\begin{equation}\label{eq:lambdai}
\lambda_i^\tau \coloneqq - \sum_{k \not=i } D_k \cdot D_i -\bm{D}_\sharp \cdot D_i,
\end{equation}
where $ \bullet \cdot \bullet $ denotes the intersection multiplicity at $0\in \C^2$ and $D_i$ are the irreducible components of $\bm{D}$. The numbers $\lambda_i^\tau$ are the first of the two ingredients we need to compute the vertical indices (they are denoted as $\lambda_i$ in \cite{Nemethi2018}).

Now, fix one component $f(D_i)$ of $f(\bm{D})$ and one parametrization (normalization)
\begin{equation*}
\begin{aligned}
\gamma: (\C, 0) &\to \big(f(D_i),0\big)\subset (\C^3, 0)\\
t&\mapsto \gamma(t)
\end{aligned}.
\end{equation*}
According to \cite[Section 4.2]{Nemethi2018}, we have a splitting $g=g_1\cdot g_2$ of $\im(f)$ along $f(D_i)$, up to permutation for the twisted components. Furthermore,
\begin{equation*}
p_1(H)\big(\gamma(t)\big)=\beta_1(t)p_1(g_1)\big(\gamma(t)\big)+\beta_2(t)p_1(g_2)\big(\gamma(t)\big) 
\end{equation*}
holds at every point $\gamma(t) \in f(D_i)\setminus\left\{0\right\}$, where $\beta_1$ and $\beta_2$ are some coefficient germs and $p_1(\bullet)$ is the first order Taylor polynomial (see further details in \cite[Definition 4.6]{Nemethi2018}). The product $ \beta_1 \beta_2 $ is a well-defined meromorphic germ $( \C, 0) \to ( \C, 0)$.  Then, we define
\begin{equation}\label{eq:vfraki}
\mathfrak{v}_{f\{i\}}^\tau\coloneqq o_{-} (\beta_1 \beta_2),
\end{equation}
where $o_{-}(\bullet)$ denotes the order (i.e., the smallest negative power) of the Laurent series ($\mathfrak{v}_{f\{i\}}^\tau$ are denoted as $\mathfrak{v}_j$ in \cite{Nemethi2018}). This is the last object we need to compute the vertical indices, hence, to  determine the surgery in \cref{th:milnperpart1}, \cref{th:milnperpar2} completely.

\begin{theorem}[see {\cite[Theorem 4.9., Lemma 4.11.]{Nemethi2018}}]
\label{th:milnper par3}
The vertical indices of $ f$ along $ f(\delta_i)$ (see \cref{def:verticalindex}) can be expressed with the terms introduced in \cref{eq:vfraki,eq:lambdai} as
\begin{equation*}
\mathfrak{vi}_{f\{i\}} \coloneqq \left\{ \begin{array}{ccc}
\lambda_i^\tau  + \lambda_{i'}^\tau +  \mathfrak{v}_{f\{i\}}^\tau  & \text{ if } & f\{i\}=\left\{i,i'\right\}, \\
\lambda_i^\tau  + \mathfrak{v}_{f\{i\}}^\tau  & \text{ if } & f\{i\} = \left\{i\right\}. \\
\end{array} \right.
\end{equation*}
\end{theorem}

Note that both $ \lambda_i^\tau$ and $\mathfrak{v}_{f\{i\}}^\tau$ depend on the choice of the transverse section $\tau$. However,  $\mathfrak{vi}_{f\{i\}}$ does not depend on it, cf. \cite[Corollary 4.12]{Nemethi2018}, thus it is an invariant of $g$ and the component $f(\delta_i)$.

Furthermore, the sum of the vertical indices is determined by the following theorem. 

\begin{theorem}\label{th:sumvert}
	For a finitely determined germ $f:(\CC^2,0)\to(\CC^3,0)$,
 \begin{equation*}
 \sum_{f\{i\}} \mathfrak{vi}_{f\{i\}} = -  \sum_{i \neq k} D_i \cdot D_k - C(f) +3T(f).
 \end{equation*}
\end{theorem}

In special cases, the formula is stated as \cite[Proposition 5.1.1.]{Nemethi2018} and proved in an algebraic way. The proof of the general formula is based on the topological description of the vertical indices $ \mathfrak{vi}_{f\{i\}} $ presented in \cite[Corollary 5.1.5]{Pinter2023}, cf. also \cref{itemd:prsumm} of \cref{pr:summ} below.
\newline

\subsection{Slice singular manifolds of immersions}\label{sec:ssmoi}
\textsl{Slice surfaces} of knots are a common concept in low-dimensional topology (see, e.g., \cite{Ozsvath2003}), they are the natural generalization of \textsl{slice discs} to any genus. More precisely, a \textit{slice surface} for a knot $\kappa\subset \Sbb^3=\partial B^4$ is a smooth submanifold of $B^4$ whose boundary is $\kappa$.
Now, we comment on a generalized version of 
slice surfaces for immersions (see \cref{rm:notationssm} below):
\textsl{slice singular manifolds}. They contain a lot of information of the \textsl{regular homotopy class} of the immersion they bound. Good examples of this and similar ideas are the works of Hughes and Melvin \cite{Hughes1985}; Saeki, Sz\H{u}cs and Takase \cite{Saeki2002}; Ekholm and Sz\H{u}cs \cite{Ekholm2003,Ekholm2006}; Takase \cite{Takase2007}; Ekholm and Takase \cite{Ekholm2011}; Kinjo \cite{Kinjo2015}; and Juhász \cite{Juhasz2005}. For a summary of this topic, see \cite[Chapter 2]{Gergothesis}. In particular, we will use the formula of Ekholm--Sz\H{u}cs--Takase--Saeki (ESzTS formula),
see \cref{eq:ib}, given in \cite[Definition 7 and Theorem 5]{Saeki2002}.
\newline

Recall that two immersions are called \textit{regular homotopic} if they are homotopic through immersions. We consider immersions of a closed oriented 3-manifold $N^3$ to $\Sbb^5$ with trivial normal bundle.


\begin{definition}\label{def:slice manifold}
A \textit{slice singular manifold} for an immersion $n: N^3 \looparrowright \Sbb^5$ is a stable smooth map $m:M^4\rightarrow B^6$  of an oriented 4-manifold with boundary $M^4$ to the closed ball $B^6$ such that
\begin{enumerate}[label=(\roman*),font=\itshape]
   \item  $\partial M^4$ is diffeomorphic  to $N^3$,
    \item $m|_{\partial M^4}$ is a stable immersion of $\partial M^4$ to $\Sbb^5$ which is regular homotopic to $n$ in $\Sbb^5$,
    \item $m^{-1}(\Sbb^5)=\partial M^4$, and
    \item $m$ is non-singular near the boundary.
\end{enumerate}

\end{definition}

\begin{remark}\label{rm:notationssm}
As far as we know, it is the first time the terminology of \textsl{slice singular manifold} is used. In many texts, one can find the terminology of \textit{singular Seifert surface} (despite it being a \textsl{manifold} instead of a \textsl{surface}) for manifolds that also bound an immersion but have ambient space the same $\Sbb^5$ instead of the ball, in the same way the classical Seifert surface lives in $\Sbb^3$ for a classical link. However, the term \textsl{singular Seifert surface} is used in \cite{Nemethi2018,Pinter2023a} instead of \textsl{slice singular manifolds}.
\end{remark}

A smooth stable map $f$ of a 4-manifold to the 6-space has very restricted type of singularities, namely:
\begin{itemize}
	\item regular simple points,
   \item surfaces of transverse double values with regular branches,
   \item isolated transverse triple values of regular branches, and
   \item curves of generalized Whitney umbrella points, which have local form
   \begin{equation}\label{eq:localWU}
   \begin{aligned}
     ( \mathbb{R}^{2} \times \R, 0) &\to (\mathbb{R}^{2} \times \mathbb{R}^{2} \times \R, 0) \\
              ( \underline{x}, y)&\mapsto( \underline{x}, y \underline{x}, y^2).
   \end{aligned}
   \end{equation}
\end{itemize} 

The ESzTS formula contains several terms: the signature of a 4-manifold $\sigma$, an invariant of a 3-manifold $\alpha$, the algebraic number of triple values of a stable map $t$, and two invariants of a stable map given as linking numbers $\ell$ and $L$. We introduce them now.

The signature of a 4-manifold $M^4$, $\sigma(M^4)$, is just the signature of the intersection form
\begin{equation}\label{eq:intersectionform}
\begin{tikzcd}[ampersand replacement = \&]
H^2(M^4,\partial M^4;\ZZ)\times H^2(M^4,\partial M^4;\ZZ) \arrow[r, "Q"] \& \ZZ                                                         \\[-0.7cm]
\hspace{2cm}(\alpha,\beta) \hspace{1.7cm}\arrow[r, maps to]                                  \& \big\langle \alpha\smile\beta,\left[M^4\right]\big\rangle.
\end{tikzcd}
\end{equation}

The closure of the set of double values of a slice singular manifold $m:M^4\rightarrow B^6$, which we denote by $m(D)$, is an immersed manifold with boundary. As an immersed manifold, $m(D)$ also has triple self-intersection points at the triple values of $m$. Each triple value is endowed with a sign depending on whether the product of the orientation of the branches of $m(D)$ agrees with the orientation of $B^6$ (positive) or not (negative). We define the \textit{algebraic number of triple values} of $m$ as the signed sum of triple values, and we denote it by $t(m)$.
\newline

\begin{figure}
	\centering
		\includegraphics[width=1\textwidth]{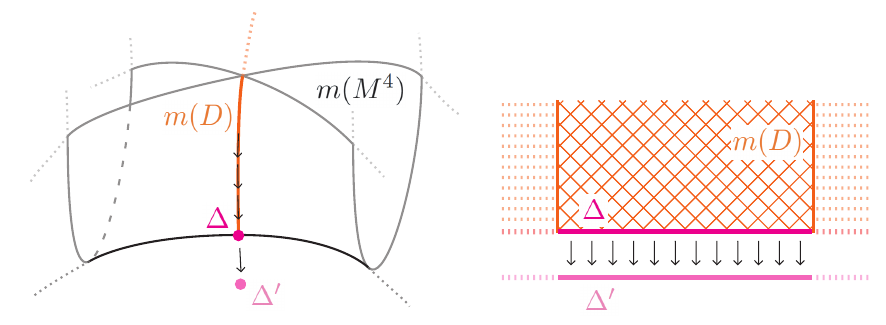}
	\caption{Two representations of the pushing out \textcolor[rgb]{0.95,0.53,0.71}{$\Delta'$} of the generalized cross-caps \textcolor[rgb]{0.93,0,0.55}{$\Delta$} from the set of double values \textcolor[rgb]{0.95,0.44,0.13}{$m(D)$} of $m(M^4)$ (cf. $\ell(m)$ in \cref{def:ell}).}
	\label{fig:l-invariant}
\end{figure}

Furthermore, the boundary $\partial m(D)$ has two parts: $m(D)\cap \Sbb^5$ and $\partial m(D) \cap \mathring{B}^6$. It is obvious that $m(D)\cap \Sbb^5$ is the set of double values of $m|_{\partial M^4}$ and, away from $\Sbb^5$, we have the set of non-immersive values of $m$, $\Delta$, given by generalized Whitney umbrellas (see in \cref{eq:localWU} that the Whitney umbrellas are one of the boundary components of the double values). Let $\Delta'$ be a copy of $\Delta$ shifted slightly along the outward normal vector field of $m(D)$. Then, $\Delta'\cap m(M^4)=\varnothing$. We define their linking number (see \cref{fig:l-invariant})
  \begin{equation}\label{def:ell}
      \ell(m)\coloneqq\lk_{(B^6, \Sbb^5)} \big(\Delta', m(M^4)\big)
  \end{equation}
  in $B^6$, relative to the boundary (see \cref{re:linkingdef}).

  The two remaining terms, $L$ and $\alpha$, of the ESzTS formula, see \cref{eq:ib} below, are not used in our argument. However, we define them to state the formula and to summarize some closely related results.  

  To define  $L$, consider a \textsl{stable} immersion of a 3-manifold $N^3$, $s: N^3 \looparrowright \Sbb^5$, with trivial normal bundle. The double values of $s$ form a manifold $\delta \subset \Sbb^5$ of dimension 1.  We can assume that the branches are orthogonal, taking a regular homotopy  of $s$ through stable immersions if necessary. 
  
  Let $v$ be a nowhere-zero section of the normal bundle of $s$. At a double value $q=s(p_1)=s(p_2)$, define $w(q)=v(p_1)+v(p_2)$. It provides a vector field $w$ along $\delta$, which is not tangent to the branches. Let $\delta'$ be the shifted copy of $\delta$ slightly along $w$. Then $\delta'$ is disjoint from the image $s(N^3)$ of $s$. We define
  \begin{equation}\label{def:L}
      L(s)\coloneqq\textnormal{lk}_{\Sbb^5} \big(\delta', s(N^3)\big).
  \end{equation}
  $L(s)$ does not depend on the choice of the normal framing $v$, see \cite[Remarks 2.6. and 5.10]{Saeki2002}.
\newline

   Finally, for a 3-manifold $N^3$ we define $$\alpha (N^3) := \dim_{\Z_2} \big(\tau H_1 (N^3, \Z) \otimes \Z_2\big),$$ where $ \tau H_1 (N^3, \Z)$ is the torsion subgroup of $ H_1 (N^3, \Z) $.
\newline

    We are ready to state our main tool (see \cite[Definition 7 and Theorem 5]{Saeki2002}).

  \begin{theorem}[The Ekholm--Sz\H{u}cs--Takase--Saeki formula, ESzTS formula]\label{th:saekiszucs}
Let $n:N^3\looparrowright \Sbb^5$ be an immersion of a 3-manifold $N^3$ with trivial normal bundle and $m:M^4\to B^6$ a slice singular manifold of $n$. The value of the expression
\begin{equation}\label{eq:ib}
    i_b (n) =  \frac{3}{2} \big( \sigma (M^4) - \alpha (N^3)\big) +  \frac{1}{2} \big(3 t(m) - 3 \ell(m) -
 L(m|_{\partial M^4})\big) 
\end{equation}
is an integer, which depends only on the regular homotopy class of the immersion $n$. 
In particular, it does not depend on the choice of the slice singular manifold $m$.
  \end{theorem}

\begin{remark}\label{re:aftersaekiszucs}
The normal bundle is always trivial for embeddings $n: N^3 \hookrightarrow \Sbb^5$ of an arbitrary 3-manifold $N^3$, because of the existence of a Seifert manifold $M^4 \subset \Sbb^5$ with boundary $N^3 \subset \Sbb^5$. We will apply the theorem in this particular case, for the inclusion of the boundary of the Milnor fiber $\partial F \subset \Sbb^5$.

  The normal bundle is also trivial for sphere immersions $n: \Sbb^3 \looparrowright \Sbb^5$, since any oriented real vector bundle of rank 2 over $\Sbb^3$ is trivial.   Indeed,
    \cref{th:saekiszucs} was stated originally for  sphere immersions in \cite{Ekholm2003}. In this case $i_b(n)$ agrees with the \textsl{Smale invariant} $\Omega(n)$, which  is a \textsl{complete} regular homotopy invariant of sphere immersions (see  \cite{Smale1959} and \cite[Theorem 1.1]{Ekholm2003}, cf. \cref{pr:summ} below). 
    
    In the general case, for any 3-manifold $N^3$ instead of $\Sbb^3$, $i_b(n)$ is not a complete regular homotopy invariant, but if we also consider the so-called \textsl{Wu invariant}, they are a complete set of invariants \cite[Theorem 5]{Saeki2002}. \cref{th:saekiszucs} is generalized for immersions $n: N^3 \looparrowright \Sbb^5$ with non-trivial normal bundle in \cite{Juhasz2005}.

\end{remark}

\begin{remark}
    In \cite{Ekholm2003, Saeki2002} one can find $L$ with positive sign on the a right hand side of \cref{eq:ib}. This is because the Ekholm--Sz\H{u}cs invariant $L(s)$ of stable immersions $s: \Sbb^3 \looparrowright \Sbb^5$ was given by an essentially different definition (see ,e.g., \cite{Ekholm2001}), causing a slight confusion in the literature. The detailed analysis in \cite{Pinter2023a, Pinter2023} shows that the two definitions agree with opposite sign. However, the definition we present here is the only possible way to generalize it for stable immersions of 3-manifolds, instead of 3-spheres, hence the minus sign for $L$ to obtain the right formulation of \cref{eq:ib}. In any case, $L$ will be 0 in our computations.
\end{remark}

\subsection{Slice singular manifold of the associated immersion}\labelpar{ss:assoc} 

We return to $\Ascr$-finite holomorphic germs  $f|_{\Sbb^3}: \Sbb^3 \looparrowright \Sbb^5$, recall \cref{sec: intro germs}. Consider the stable immersion $f|_{\Sbb^3}: \Sbb^3 \looparrowright \Sbb^5$ associated to $f$,  as shown in \cref{ss:associmm}.  In \cite[Section 9]{Nemethi2015}, a slice singular manifold for $f|_{\Sbb^3}$ is created in the sense of \cref{sec:ssmoi}, by a stabilization of $f$ in the real $\mathcal{C}^{\infty}$ sense. This construction is summarized here, including the properties of this slice singular manifold appearing in \cref{th:saekiszucs}.

We start with a stable deformation  $f_s$ of $f$ in the holomorphic sense, see \cref{sec: intro germs}. This defines a map $f_s: B^4 \to B^6_{\epsilon}$, where $B^6_{\epsilon}$ is the fixed Milnor ball of $(\mathcal{S}, 0)=(\im f, 0) \in (\C^3, 0)$, and its preimage $f_s^{-1}(B^6_{\epsilon})$ is diffeomorphic to $B^4$. Moreover, $f_s^{-1}(B^6_{\epsilon})$ can be naturally identified with the ball $f^{-1}(B^6_{\epsilon})$, see \cite[Section 9.1]{Nemethi2015}.

\begin{figure}
	\centering
		\includegraphics[width=1\textwidth]{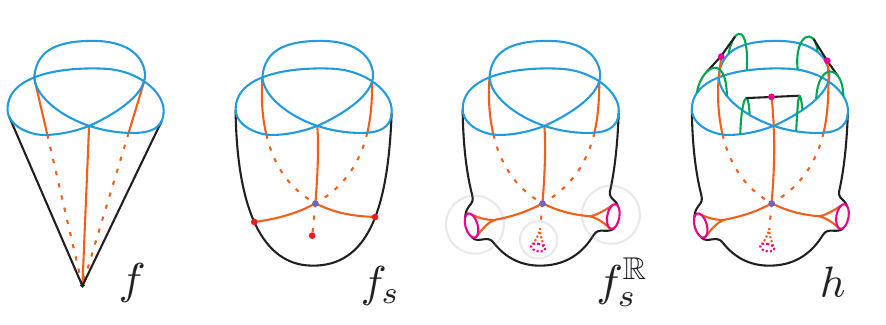}
	\caption{An $\eqA$-finite map germ $f$; a holomorphic stable deformation $f_s$ with complex \textcolor[rgb]{1,0,0}{cross-caps} and \textcolor[rgb]{0.39,0.44,0.73}{triple values}; a $\mathcal{C}^{\infty}$ stable deformation $f_s^\RR$  with curves of \textcolor[rgb]{0.93,0,0.55}{generalized real cross-caps}; and the map $h$ we construct, with new curves of \textcolor[rgb]{0.93,0,0.55}{generalized real cross-caps}. Note that $f_s^\RR$ and $h$ are slice singular manifolds of the associated immersion \textcolor[rgb]{0,0.6,0.85}{$f|_{\Sbb^3}$} and of $\partial F \subset \Sbb^5$, respectively. See \cref{ss:assoc,ss:constructionh}.}
	\label{fig:4figs}
\end{figure}

$f_s$ is stable in the complex analytic sense, but it is not stable as a real smooth map. Indeed, comparing the list of stable germs of holomorphic maps $\C^2 \to \C^3$ and $\mathcal{C}^{\infty}$ maps $\R^4 \to \R^6$ (see \cref{sec: intro germs} and \cref{sec:ssmoi}), we see that the only obstruction of $f_s$ to be stable in the smooth sense is that it has complex cross-caps. However, we can modify $f_s$ to a smooth stable map 
$$f_s^\RR:B^4\to B^6$$
in the following way (see also \cref{fig:4figs}).

We choose local complex coordinates for $f_s$ such that, around a complex cross-cap, we have the normal form $f_s(u,v)=(u,v^2,uv)$ (recall \cref{ex:cc}). We take its real $\mathcal{C}^{\infty}$ deformation, with $0<r\ll|s|$,
\begin{equation*}
(u,v)\mapsto(u,v^2+2r\overline{v},uv+r\overline{v}).
\end{equation*}
This deforms the complex cross-cap to an $\Sbb^1$-family of generalised real cross-caps (cf. \cref{eq:localWU}), as it is shown in \cite[Section 10.1]{Nemethi2015}. This implies that it is a stable map in the real $\mathcal{C}^{\infty}$ sense.
Then, since away from a complex cross-cap $f_s$ is stable in the smooth sense, we can take partitions of unity to glue the previous local $\Ccal^\infty$ modification to $f_s$, having the $\Ccal^\infty$ stable map $f_s^\RR$. 

As we mentioned in \cref{ss:associmm}, the immersion  $f_s|_{\Sbb^3}$  and $f|_{\Sbb^3}$ are regular homotopic through stable immersions. Moreover, the change from $f_s$ to $f_s^\RR$ leaves a neighbourhood of the boundary (in the source and target) unchanged. Therefore $f_s^\RR: B^4 \to B^6$ defines a slice singular manifold of the associated immersion $f|_{\Sbb^3}$.

\begin{remark}\label{re:identification}
What is more, the balls, $f^{-1}(B^6_{\epsilon})$, $f^{-1}_s(B^6_{\epsilon})$ and $(f^{\R}_s)^{-1}(B^6_{\epsilon})$ and their boundaries are canonically isotopic in $\C^2$. For the construction in \cref{s:ssmftbotmf} below, it is convenient to identify them and refer to them simply by $B^4$ (which serves as a common domain of the three maps) with boundary $\Sbb^{3}$ (the domain of the associated immersion). Moreover, the identification can be done in 
 such a way that the maps agree near the boundary $\Sbb^3$. In  particular, we have the fiber bundle structure described in \cref{lem:twORuntw} also for $f_s^\RR$.
 \end{remark}

 \begin{proposition}\label{pr:summ} 
 
 The terms of the ESzTS formula \cref{th:saekiszucs} applied to $f_s^\RR$ can be expressed in terms of $C(f)$ and $T(f)$ as follows:

  \begin{enumerate}[label=(\alph*)]     \item\label{itema:prsumm} $\Omega(f|_{\Sbb^3})=-C(f)$,
   \item\label{itemb:prsumm}   $t(f_s^\RR)=T(f)$,
   \item\label{itemc:prsumm} $\ell(f_s^\RR)=C(f)$,
   \item\label{itemd:prsumm} $L(f|_{\Sbb^3})=3T(f)-C(f)$.
 \end{enumerate}
%
%
%
 \end{proposition}

Recall that $\Omega(f|_{\Sbb^3})$ is the Smale invariant of sphere immersions, which agrees with $i_b(f|_{\Sbb^3})$. Despite this,
\cref{itema:prsumm} (known as \textsl{holomorphic Smale invariant formula}) is proved in \cite[Theorem 1.2.2]{Nemethi2015} independently of \cref{th:saekiszucs}. \cref{itemb:prsumm} follows trivially from the construction. The proof of \cref{itemc:prsumm} in \cite[Theorem 9.1.3]{Nemethi2015} is based on a local calculation for the $\mathcal{C}^{\infty}$-stabilization of the complex Whitney umbrella we have shown above. Then, \cref{itemd:prsumm} is deduced by substituting the other terms in the formula, completed with local computations on examples. An independent, direct proof for \cref{itemd:prsumm} is the main result of \cite{Pinter2023a}. 

\begin{remark}
    In \cite{Nemethi2018},  \cref{itemd:prsumm} of \cref{pr:summ} is stated with opposite sign, i.e., as `$L(f_s^\RR)=C(f)-3T(f)$'. This is because the other definition of $L$ is used, which results the same invariant with opposite sign. For details we refer to \cite[Appendix]{Pinter2023a}. Since we do not need $L$, this ambiguity is not a problem, cf. \cref{sss:linkinglL}.
\end{remark}

\section{Slice singular manifolds for the boundary of the Milnor fiber}\label{s:ssmftbotmf}

\subsection{Preliminary summary of the construction}

We follow the general idea of Andr\'as N\'emethi of using \cref{th:saekiszucs} to compute the signature of the Milnor fiber by finding a suitable slice singular manifold for the boundary. In this paper we fulfil this for images of finitely determined map germs.
\newline

Given an $\mathscr{A}$-finite holomorphic map germ $f: (\C^2, 0) \to (\C^3, 0)$, its image $(\mathcal{S},0)=\big(\im(f),0\big)\subset(\CC^3,0)$ and its Milnor fiber $F_f$, we construct a slice singular manifold $h:X\to B^6$ for the embedding $\partial F_f \subset \Sbb^5$, constructed from the germ $f$.

Since we need a stable map in the smooth sense (see \cref{def:slice manifold}), we start with $f_s^\RR$ given in \cref{ss:assoc}. Indeed, this is a slice singular manifold for the immersion $f|_{\Sbb^3}$. We extend its domain $B^4$ by taking the trace of the surgery described in \cref{th:milnperpart1} to get a manifold $X$ with $\partial X \cong \partial F_f$.
We also extend $f_s^\RR$ to a map $h:X\to B^6$, with $h|_{\partial X}$ regular homotopic to $\partial F_f \hookrightarrow \Sbb^5$ in $\Sbb^{5}$. Even more is true for our particular construction: $h|_{\partial X}$ is also an embedding and the two embeddings are isotopic.

The invariant given in \cref{th:saekiszucs} can be computed by using two different slice singular manifolds: the embedding $F_f 
 \hookrightarrow B^6$, and the map $h: X \to B^6$. Comparing the two expressions, without any further analysis of the slice singular manifold $h$, we have
 \begin{equation*}
     \sigma(F_f)=\sigma(X) +t(h)-\ell(h).
 \end{equation*}
 
 For our particular construction, we show that new triple values are not created during the extension of $f_s^{\R}$ to $h$, and, although new singular points (generalised real cross-caps) are created, they do not change $\ell$. By these properties and \cref{pr:summ}, we get the values $t(h)=t(f_s^{\R})=T(f)$ and $\ell(h)=\ell(f_s^{\R})=C(f)$, proving \cref{th:main}.

\subsection{Construction of the 4-manifold}\label{ss:4mfld}

The 4-manifold $X$ is defined as the trace of the surgery in \cref{eq:rag} used for the construction of $\partial F$. To describe it, we have to introduce the 4-manifold $Z$, the solid version of the 3-manifold $Z_\partial$ from \cref{eq:Y}, as
\begin{equation}\label{eq:3Y}
Z\coloneqq \frac{\left[-\pi,\pi\right]\times B^2\times\left[-1,1\right]}{(-\pi,z,t)\sim(\pi,\overline{z},-t)}.
\end{equation}
Observe the following properties of $Z$. 

\begin{proposition}\label{prop:Z}
 For $Z$ defined as in \cref{eq:3Y},
\begin{enumerate}[label=(\roman*)]
\item\label{it:Z1} $Z $ is diffeomorphic to $\Sbb^1 \times B^3$.
\item\label{it:Z2} The boundary $\partial Z \cong \Sbb^1 \times \Sbb^2$ decomposes as
\begin{equation*}
\partial Z=Z_\partial \cup_{\partial} (\Sbb^1 \times B^2),
\end{equation*}
with gluing along the torus boundaries $\partial Z_\partial \cong \Sbb^1 \times \Sbb^1$.
\item\label{it:Z3} The centers of the discs $B^2$ form a Möbius band, that is, 
\begin{equation}\label{eq:mobiusband}
\mathcal{M}=\frac{\left[-\pi,\pi\right] \times \{0\} \times \left[-1,1\right]}
{(-\pi,0,t)\sim(\pi,0,-t)} \subset Z.
\end{equation}
\end{enumerate}
\end{proposition}

\begin{proof} The projection $Z\to \Sbb^1$ to the first component defines a locally trivial fibration with fibers $B^2 \times \left[-1,1\right] \cong B^3$. The geometric monodromy is $(z, t) \mapsto (\bar{z}, -t)$ and it is isotopic to the identity of $B^3$, hence the fibration is trivial. This proves \cref{it:Z1}. 

The projection restricted to the boundary defines a locally trivial fibration $\partial Z \to \Sbb^1$ with fibers 
\begin{equation*}
\partial \big(B^2 \times \left[-1,1\right]\big)=\big(\Sbb^1 \times \left[-1,1\right]\big) \cup_{\partial} \big(B^2 \times \left\{-1,1\right\}\big).
\end{equation*}
The total space corresponding to $\Sbb^1 \times \left[-1,1\right]$ is equal to $Z_\partial$ by definition (recall \cref{eq:Y}). The other part is diffeomorphic to a solid torus $\Sbb^1 \times B^2$:
\begin{equation*}
\frac{\left[-\pi,\pi\right]\times B^2\times\left\{-1,1\right\}}{(-\pi,z,t)\sim(\pi,\overline{z},-t)}.
\end{equation*}
This proves \cref{it:Z2}. \cref{it:Z3} is trivial.
\end{proof}


Then, recovering the construction from \cref{ss:surgerypartialf}, define the pieces 
\begin{equation}\label{eq:Ziinsecondform}
    Z_{f\{i\}} :=
    \left\{
\begin{array}{cc}
\Sbb^1 \times B^2 \times [-1,1] &  \textnormal{ if } f\{i\}=\{i, i'\} \\

Z & \textnormal{ if } f\{i\}=\{i\},
\end{array}
    \right.
\end{equation}
recall \cref{not:fi,eq:3Y}.
Then, 
$X$ is defined as
\begin{equation}\label{eq:def of X}
X\coloneqq\left[ B^4 \sqcup  \bigsqcup_{f\{i\}} Z_{f\{i\}}  \right]  \Bigg / \Phi
\end{equation}
 with a gluing $\Phi$ that extends $\Phi_\partial$ of the construction of $\partial F$ in \cref{th:milnperpart1} in the following sense (see also \cref{fig:X2d,fig:X3d}). $\Phi$ is induced by the collection of diffeomorphisms 
\[ \Phi^{f\{i\}}:\begin{cases}
 N_i \sqcup  N_{i'} \to  \Sbb^1 \times B^2 \times \{-1,1 \} \subset \partial (\Sbb^1 \times B^2 \times [-1,1 ]) & \text{, }f\{i\}=\left\{i,i'\right\}\\
  N_i \to  \Sbb^1\times B^2\subset\partial Z & \text{, }f\{i\}=\left\{i\right\},
\end{cases} \]
where $\Phi^{f\{i\}}$ are defined in the same way as in \cref{th:milnperpar2}. Namely, the mapping gluing $ \partial N_i $ and $ \partial N_{i'} $ (or $  \partial Y_{f\{i\}}=\partial Z_{\partial} \cong \Sbb^1\times \Sbb^1 $) given by the matrix	\begin{equation}\label{eq:gluingmatrix}
	\left( \begin{array}{cc}
	-1 & \mathfrak{vi}_{f\{i\}} \\
	0 & 1 \\
	\end{array} \right)
	\end{equation}              
   can be extended to a mapping gluing $N_i$ and $N_{i'}$ (or $\Sbb^1\times B^2\subset\partial Z$), since we map meridian to meridian. This extension induces $\Phi^{f\{i\}}$.

\begin{remark}\label{rem:othergluing}
For each untwisted component, the gluing $\Phi^{f\{i\}}$ can be replaced by an identification $N_i \to -N_{i'}$ of the two solid tori. It gives $X$ in the alternative form
\begin{equation}\label{eq:def of X2}
X\cong \left[ B^4 \sqcup  \bigsqcup_{f\{i\}=\{i\}} Z_{f\{i\}}  \right]  \Bigg / \Phi,
\end{equation}
where $Z_{f\{i\}} \cong Z$ for all twisted components and $\Phi^{f\{i\}}$ as described above, but with  $\Phi^{f\{i\}}: N_i \to -N_{i'}$ given by the gluing matrix \cref{eq:gluingmatrix} for untwisted components.

\end{remark}

\begin{proposition}\label{prop:xcongf}
From the construction above, we have
    $\partial X \cong \partial F_f$.
\end{proposition}

\subsection{Signature of $X$}\label{ss:signature} The signature of $X$ is of central importance for our purpose. To compute it, we consider $X$ in the form of \cref{eq:def of X2}.

Observe that $X$ is homotopy equivalent with its 2-skeleton. That is, the wedge of surfaces $\mathcal{C}_{f\{i\}}$ constructed as follows:
\begin{itemize}
	\item For an untwisted component $f\{i\}=\{i,i'\}$,
\begin{equation}\label{eq:Cii'}
  \mathcal{C}_{f\{i\}} \coloneqq  D_i \cup_{\Phi} D_{i'}.
\end{equation}
It is homeomorphic to $\Sbb^2$ with two opposite points glued together.
\item For a twisted component $f\{i\}=\{i\}$ (recall \cref{eq:mobiusband}),
\begin{equation}\label{eq:Ci}
  \mathcal{C}_{f\{i\}} \coloneqq  D_i  \cup_\Phi \mathcal{M} .
\end{equation}
It is homeomorphic to a real projective plane $\R P^2$ since it is a disc and a Möbius band glued along their boundaries. 
\end{itemize}
 Therefore,
\begin{equation}\label{eq:HkX}
\widetilde{H}_k(X,\ZZ)\cong
\begin{cases}
 \bigoplus_{f\{i\}=\{i,i'\}} \ZZ \oplus \bigoplus_{f\{i\}=\{i\}}\nicefrac{\ZZ}{2\ZZ} & k=1\\
 \bigoplus_{f\{i\}=\{i,i'\}}\ZZ & k=2\\
 0 & \textnormal{otherwise}.
\end{cases}
\end{equation}
%

\begin{figure}
	\centering
		\includegraphics[width=1.0\textwidth]{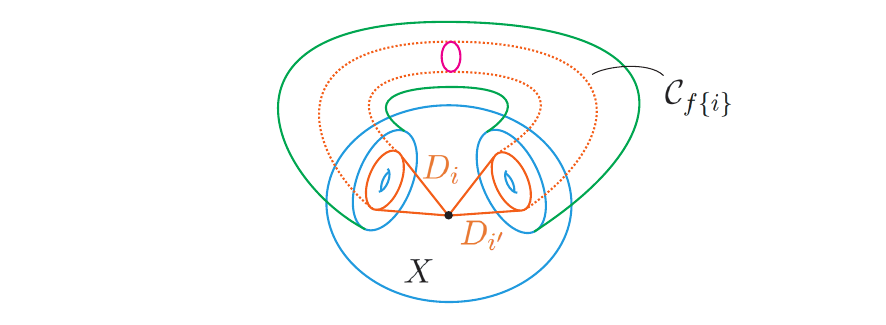}
	\caption{The homological cycle $\mathcal{C}_{f\{i\}}$ in $X$ corresponding to a couple of untwisted components \textcolor[rgb]{0.95,0.44,0.13}{$D_i$} and \textcolor[rgb]{0.95,0.44,0.13}{$D_{i'}$}, cf. \cref{eq:Cii'}. The \textcolor[rgb]{0.95,0.44,0.13}{cylinder (dashed)} connecting them can be  simplified to an identification of the boundaries (cf. \cref{eq:def of X,eq:def of X2}).}
	\label{fig:X3d}
\end{figure}

Finally, we describe the (homological) intersection form of $X$ (cf. \cref{eq:intersectionform,fig:X3d}).

\begin{proposition}\label{prop:intersection matrix}
    For different untwisted components $f\{i\}=\{i,i'\} \neq f\{j\}=\{j,j'\}$, the intersection number of the two cycles $[\mathcal{C}_{f\{i\}}]$ and $[\mathcal{C}_{f\{j\}}]$ is
\begin{equation}\label{eq:diffint}
    [\mathcal{C}_{f\{i\}}] \cdot [\mathcal{C}_{f\{j\}}]=D_i\cdot D_j+D_i\cdot D_{j'}+D_{i'}\cdot D_j+D_{i'}\cdot D_{j'},
    \end{equation}
using algebraic intersection multiplicities.

    The self-intersection of a cycle $\mathcal{C}_{f\{i\}}$ is 
    \begin{equation}\label{eq:diffint2}
     [\mathcal{C}_{f\{i\}}] \cdot [\mathcal{C}_{f\{i\}}]=2 D_i \cdot D_{i'}+ \mathfrak{vi}_{f\{i\}}.
     \end{equation}
\end{proposition}
\begin{proof}
By \cref{eq:Cii'}, we have
   $ \mathcal{C}_{f\{i\}} =  D_i \cup_{\Phi} D_{i'}
   $.
To prove \cref{eq:diffint} observe that, in general, the surfaces $\mathcal{C}_{f\{i\}}$ and $\mathcal{C}_{f\{j\}}$ have non-transverse intersection at the origin, and nowhere else. To determine their intersection number, we change $\mathcal{C}_{f\{i\}}$ to a homologous cycle. We can replace $D_i$ and $ D_{i'}$ 
by any slice surfaces of $\delta_i$ and $\delta_{i'}$ (see \cref{sec:ssmoi}), say $\Lcal_i$ and $\Lcal_{i'}$. For example, one can take the Milnor fibers $F_i,F_{i'}$ of $D_i,D_{i'}$ and glue the trace of the isotopies that take $\partial F_i,\partial F_{i'}$ to $\delta_i,\delta_{i'}$; since these last two cylinders lie in $\Sbb^3$, we can push them to the interior of $B^4$ and then we have slice surfaces as intended. Then, $\Lcal_i \cup_{\Phi} \Lcal_{i'}$ is homologous to $\mathcal{C}_{f\{i\}}$.

Since we can modify $\Lcal_i \cup_{\Phi} \Lcal_{i'}$ to have transverse intersection with $\mathcal{C}_{f\{j\}}=D_{j} \cup_{\Phi} D_{j'}$, the intersection number $  [\mathcal{C}_{f\{i\}}] \cdot [\mathcal{C}_{f\{j\}}]$ is equal to the algebraic number of transverse intersection points  of $\Lcal_i \cup_{\Phi} \Lcal_{i'}$ and $\mathcal{C}_{f\{j\}}=D_{j} \cup_{\Phi} D_{j'}$, which agrees with the right-hand side of \cref{eq:diffint}. Notice that it is also equal to the linking number of the links $\delta_i \sqcup \delta_{i'}$ and $\delta_j \sqcup \delta_{j'}$ in $\Sbb^3$.

We can do a similar trick to prove \cref{eq:diffint2} for the self-intersection 
    $$[\mathcal{C}_{f\{i\}}] \cdot [\mathcal{C}_{f\{i\}}].$$
We start by changing $D_i$ to its Milnor fiber $F_i$. Its boundary is isotopic to the fixed topological longitude of $\partial N_i$ (cf. \cref{th:milnperpart1,th:milnperpar2}), by an isotopy in $N_i \setminus \delta_i$.
This follows from the fact that the linking number of $\partial F_i $ and $\delta_i$ in $\Sbb^3$ is 0. Now, let $\mathcal{F}_i$ be the 2-chain obtained by gluing together $F_i$ and the trace of this isotopy along $\partial F_i$.

   By the construction of \cref{th:milnperpart1,th:milnperpar2}, in particular \cref{eq:milnpermatrix}, the topological longitude $\partial \mathcal{F}_i \subset \partial N_i$ glues to a curve $\overline{\delta}_{i'} \subset \partial N_{i'}$ with linking number 
   \begin{equation}\label{eq:linkdeltadeltabar}
       \textnormal{lk}_{\Sbb^3} (\overline{\delta}_{i'}, \delta_{i'})=\mathfrak{vi}_{f\{i\}}.
   \end{equation}

   Let $\Lcal_{i'} \in B^4$ be an arbitrary slice surface of $\overline{\delta}_{i'}$. Then, $\mathcal{F}_i \cup_{\Phi} \mathcal{L}_{i'} $ is a 2-chain in $X$ homologous with $\mathcal{C}_{f\{i\}}$. Hence,
    $[\mathcal{C}_{f\{i\}}] \cdot [\mathcal{C}_{f\{i\}}]$ is equal to the algebraic number of transverse intersection points of $\mathcal{C}_{f\{i\}}$ and $\mathcal{F}_i \cup_{\Phi} \mathcal{L}_{i'} $, which is determined by the following properties:
   \begin{itemize}
      \item The intersection of $\mathcal{F}_i$ with $D_i$ is empty, since $F_i$ is a Milnor fiber of $D_i$, and the isotopy does not cross $\delta_i$.
      \item The intersection of $\mathcal{F}_i$ with $D_{i'}$ is $D_i\cdot D_{i'}$.
      \item The intersection of $\Lcal_{i'}$ with $D_{i}$ is $D_i\cdot D_{i'}$, that is, the linking number of $\delta_i $ and $\delta_{i'}$.
      \item The intersection of $\Lcal_{i'}$ with $D_{i'}$ is equal to the linking number between $\bar{\delta}_{i'}$ and $\delta_{i'}$, which equals $\mathfrak{vi}_{f\{i\}}$ (by \cref{eq:linkdeltadeltabar}).
   \end{itemize}
   The sum of all these terms equals the right-hand side of \cref{eq:diffint2}, and the result follows.
\end{proof}

\subsection{Construction of the map}\label{ss:constructionh} 

\begin{figure}
	\centering
		\includegraphics[width=1.0\textwidth]{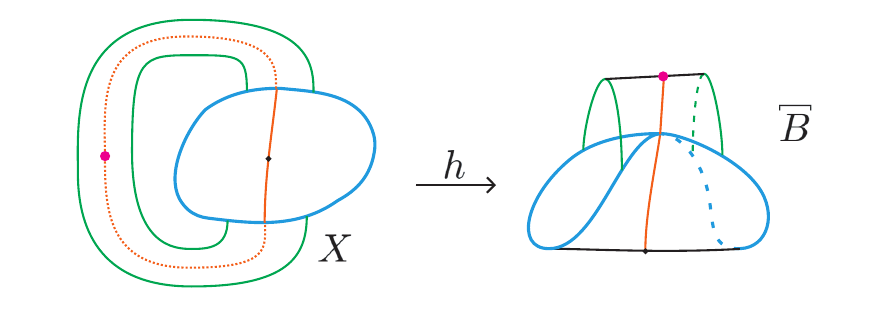}
	\caption{The map $h$ from $X$ to $\shat{B}$.}
	\label{fig:X2d}
\end{figure}

We consider $X$ defined as \cref{eq:def of X}, i.e., $B^4$ with some glued pieces. The disjoint union of these pieces is called the exterior part of $X$, $X_{\textnormal{ext}}\coloneqq \textnormal{cl}(X\setminus B^4)$. We extend $f_s^\mathbb{R}$ to a map $h$ defined on $X$ by creating one dimensional families of generalized real Whitney umbrellas on $X_{\textnormal{ext}}$ (see \cref{eq:localWU,fig:X2d}). More precisely, the locus of the generalized Whitney umbrella points will be the \textsl{middle circles} of the glued pieces, 
\begin{equation}\label{eq:whuntw}
    \Sbb^1 \times \{0 \} \times \{0\} \subset \Sbb^1 \times B^2 \times [-1, 1]
\end{equation}
for untwisted components and the middle circle of the Möbius band \cref{eq:mobiusband}
\begin{equation}\label{eq:whtw}
    \Sbb^1 \times \{0 \} \times \{0 \} \subset \mathcal{M} \subset Z
\end{equation}
 for twisted components.

The target space of $h$ will be an extension $\shat{B}$ of $B^6$, still diffeomorphic to $B^6$, such that the exterior part $X_{\textnormal{ext}}$ is mapped to \textit{the exterior part of }$\shat{B}$ 
$$
   \shat{B}_{\textnormal{ext}} \coloneqq \textnormal{cl}\big(\shat{B} \setminus B^6_\epsilon\big),$$
providing a map
$${h}:X\to\shat{B}\cong B^6 $$
whose restriction to $B^4$ is $h|_{B^4}=f^{\R}_s$ (see \cref{fig:4figs,fig:X2d}).

\subsubsection{The target space of $h$}\label{sss:targetspace} 

Consider a component $f(\delta_i)$ of $f(\bm{\delta})$, either twisted or untwisted.  In order to define the target space $\shat{B}$, consider the tubular neighbourhood $T_{f\{i\}}$ of $f(\delta_i)$ as in \cref{lem:twORuntw}. One can take a diffeomorphism $\phi^{f\{i\}}:T_{f\{i\}}\rightarrow\Sbb^1 \times B^4$ such that, by \cref{lem:twORuntw}, the image of $K_f\cap T_{f\{i\}}$ by this diffeomorphism $\phi^{f\{i\}}$  in $\Sbb^1 \times B^4$ is a locally trivial fiber bundle with fiber an $A_1$ singularity $\left\{\alpha\beta=0\right\}$, with the coordinates $(\alpha,\beta)$ of the unit ball $B^4\subset \CC^2$. Indeed, for untwisted components this is immediate. For twisted components $\phi^{f\{i\}}$ identifies $T_{f\{i\}}$ with a $B^4$-bundle over $\Sbb^1$ that carries the $A_1$-bundle of \cref{lem:twORuntw}. Hence, its geometric monodromy interchanges the coordinates $\alpha$ and $\beta$ of $B^4$. That is,
\begin{equation*}
    \frac{[-\pi, \pi] \times B^4 }
    {(-\pi;\: \alpha,\beta) \sim (\pi;\: \beta,\alpha)}.
\end{equation*}
However this bundle is obviously trivial, since its geometric monodromy $(\alpha,\beta) \mapsto (\beta,\alpha)$ is isotopic to the identity of $B^4$. Hence, for both untwisted and twisted components we have a diffeomorphism $T_{f\{i\}} \cong \Sbb^1 \times B^4$. Notice that we still have freedom to define $\phi^{f\{i\}}$ with the same properties, we will fix this freedom in \cref{lem:h01}

Using this diffeomorphism, define $\varphi$ as the collection of gluings $\varphi^{f\{i\}}$:
\begin{equation}\label{eq:phiT1}
\begin{tikzcd}
 T_{f\{i\}}\arrow[r,"\phi^{f\{i\}}"]\arrow[rr, "\varphi^{f\{i\}}",
rounded corners,
to path={ -- ([yshift=5ex]\tikztostart.center) --  ([yshift=5ex]\tikztotarget.center) \tikztonodes -- (\tikztotarget.north)}
] &\Sbb^1\times B^4 \times \{-1\}\arrow[r,"i",hook]&\Sbb^1\times B^4 \times \left[-1,1\right]
\end{tikzcd},\end{equation}
This defines the target of our final map $h$,
\begin{equation}\label{eq:defbshat}
\shat{B}\coloneqq \left[ B_\epsilon^6 \sqcup \bigsqcup_{f\{i\}} \Sbb^1\times B^4 \times \left[-1,1\right]\right]\Bigg / \varphi.
\end{equation}

As it was the case for $X$, $\shat{B}$ can be slightly deformed such that it has smooth boundary, being diffeomorphic to $B^6$.

\subsubsection{The map on the exterior part} 
\begin{figure}
	\centering
		\includegraphics[width=1.05\textwidth]{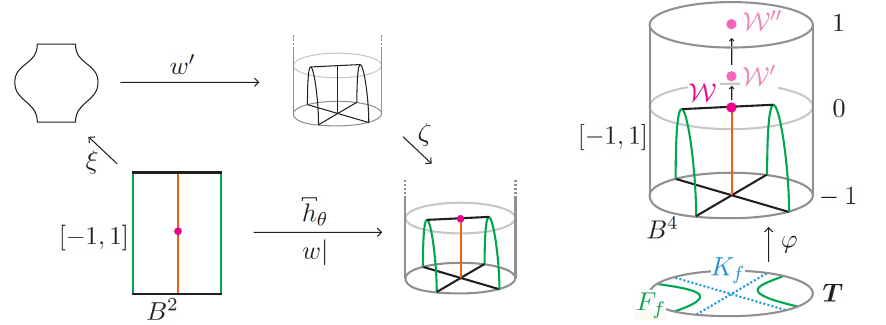}
	\caption{Left: The map $\shat{h}_\theta$, as the restriction $w|$  of a map left-right equivalent to the generalized cross-cap $w'$ by $\xi,\zeta$ (see \cref{eq:commutativeww'}). Right: The gluing $\varphi$ of \cref{eq:phiT1,eq:defbshat}. The boundary parts described in \cref{it1:propertiesofshath,it3:propertiesofshath} of \cref{pr:propertiesofshath}, $\textcolor[rgb]{0,0.53,0.22}{F_f}\cap\bm{T}$, and $\textcolor[rgb]{0.11,0.5,0.69}{K_f} \cap \bm{T}$ are also depicted (cf. \cref{re:carriesthechange,sss:themapontheboundary}). The pushing out \textcolor[rgb]{0.95,0.53,0.71}{$\mathcal{W}'$} and \textcolor[rgb]{0.95,0.53,0.71}{$\mathcal{W}''$} of the cross-cap points \textcolor[rgb]{0.93,0,0.55}{$\mathcal{W}$} of the proof of \cref{th:lext}. }
	\label{fig:cylinder}
\end{figure}

In order to define $h$, for every component $f\{i\}$ we need some pieces,
\begin{equation*}
\shat{h}|: \begin{cases}
     \Sbb^1\times B^2\times \left[-1,1\right] \longrightarrow \Sbb^1\times B^4 \times \left[-1,1\right] & f\{i\}=\{ i,i'\}\\
     Z\longrightarrow \Sbb^1\times B^4 \times \left[-1,1\right] & f\{i\}=\{i\}.
\end{cases}
\end{equation*}
The collection of these pieces defines the exterior map $X_{\textnormal{ext}} \to \shat{B}_{\textnormal{ext}}$ that will be glued to $f_s^\RR$ in a convenient way. Moreover, each piece will be trivial over $\Sbb^1$, i.e., for every $\theta\in \Sbb^1$ the map $\shat{h}_\theta\coloneqq \shat{h}|(\theta,\bullet,\bullet)$ is the same. Hence, it is enough to define $\shat{h}_\theta$. Actually, $\shat{h}_\theta$ will be left-right equivalent to the generalized cross-cap (cf. \cref{eq:localWU,fig:cylinder})
\begin{equation*}
\begin{aligned}
w':\mathbb{C}\times \mathbb{R}&\longrightarrow \mathbb{C}^2\times \mathbb{R}\\
	(z;\:t)&\longmapsto \big(z, zt;\:-t^2\big)
\end{aligned}.
\end{equation*}
However we have to modify this normal form to obtain a map $\shat{h}_\theta$ with the desired properties (collected below). First, we take a map $w$ left-right equivalent to $w'$ as follows:
\begin{equation}\label{eq:commutativeww'}\begin{tikzcd}
\CC\times\RR \arrow[r, "w" ]& \CC^2\times\RR\arrow[d, " \zeta","\sim"'  {anchor= south, rotate=90, inner sep=.3mm}] \\
\CC\times\RR \arrow[u, " \xi ","\sim"  {anchor= north, rotate=90, inner sep=.6mm}]\arrow[r, "w'" ]& \CC^2\times\RR
\end{tikzcd},\end{equation}
where we take the diffeomorphisms
\begin{equation*}
\begin{aligned}
    \xi(z; t)&=\left( \frac{z}{\sqrt{2+2t^2}}; \ t \right), \text{ and}\\
    \zeta(\alpha, \beta; s)&=\left(\alpha+\beta, \ \overline{\alpha}-\overline{\beta}; \ s
    \right)
   \end{aligned}
\end{equation*}
Indeed, $\zeta$ maps $\{(\alpha+ \beta)(\alpha-\beta)=0\}\times\left\{s\right\}$ to $\{\alpha \beta=0\}\times\left\{s\right\}$ and it reverses the orientation on one of the two discs, this is the reason to take its composition (cf. \cref{it3:propertiesofshath}, \cref{it3:propertiesofshath ori} of \cref{pr:propertiesofshath}). To be more precise, the composition is
\begin{equation}\label{eq:w}
\begin{aligned}
w:\mathbb{C}\times \mathbb{R}&\longrightarrow \mathbb{C}^2\times \mathbb{R}\\
	(z;\:t)&\longmapsto \left(
 \frac{z(1+t)}{\sqrt{2+2t^2}}, \ 
 \frac{\overline{z}(1-t)}{\sqrt{2+2t^2}}
 ;\:-t^2 \right).
\end{aligned}
\end{equation}

Now, we take the restriction of $w$ to $B^2 \times [-1, 1]$. Indeed, the reason we compose by $\xi$ is to have a \textsl{good enough} restriction, as we show now.

\begin{proposition}\label{pr:preimage}
    The preimage $w^{-1}(B^4 \times [-1, 1])$ is equal to $B^2 \times [-1, 1]$, where $B^2 \subset \C \cong \R^2$ and $B^4 \subset \C^2 \cong \R^4$ are the unit balls.
\end{proposition}
\begin{proof}
Since
   \begin{equation*}
      |(\alpha, \beta)|^2= 
\left| \frac{z(1+t)}{\sqrt{2+2t^2}} \right|^2 + 
\left|
\frac{\overline{z}(1-t)}{\sqrt{2+2t^2}}
\right|^2 
        =|z|^2,
   \end{equation*}
if $t\in[-1, 1]$, the pre-image of the unit ball $B^4 \times \{-t^2\}$ is the two copies of the unit disc $B^2 \times \{\pm t\}$.
\end{proof}

Then, we define $\shat{h}_{\theta}$ as the restriction of $w$ to $B^2 \times [-1, 1]$,
\begin{equation}\label{eq:defhtheta}
    \shat{h}_{\theta}\coloneqq w|_{B^2 \times [-1, 1]}: B^2 \times [-1, 1] \to B^4 \times [-1, 1].
\end{equation}

We summarise the properties of $\shat{h}_{\theta}$, by construction  (see also \cref{fig:cylinder}).

\begin{proposition}\label{pr:propertiesofshath}
The map $\shat{h}_{\theta}$ has the following properties:
\begin{enumerate}[label=(\roman*)]
    \item\label{it1:propertiesofshath} $\shat{h}_{\theta}$ gives an embedding of $\Sbb^1\times\left[-1,1\right]\subset\partial\big(B^2\times \left[-1,1\right]\big)$ to   $\Sbb^3\times\left[-1,1\right]\subset\partial (B^4\times\left[-1,1\right])$.
    \item\label{it3:propertiesofshath} $\shat{h}_{\theta}$ maps the two discs $B^2\times \left\{ t,-t\right\}$, $t \in\left]0,1\right]$, to two transverse discs intersecting at $(0;\: -t^2)\in B^4\times\left[-1,1\right]$. In particular, for $t=1$, their images are $\{\beta=0\}$ and $\{\alpha=0\}$ in $B^4\times\left\{-1\right\} \subset\partial (B^4\times\left[-1,1\right])$. 

\item\label{it3:propertiesofshath ori} The restriction of $\shat{h}_{\theta}$ induces an orientation-preserving map $B^2 \times \{1\} \to \{\beta =0\}$, and it reverses the orientation for $B^2 \times \{-1\} \to \{\alpha =0\}$, where each disc is considered with its complex orientation.
    \item\label{it56:propertiesofshath}  $\shat{h}_{\theta}$ has a unique non-immersive point at $(0;\: 0)\in B^2\times\left[-1,1\right]$, which is mapped to $(0;\: 0)\in B^4\times\left[-1,1\right]$, and it is generic (stable) in the $\mathcal{C}^{\infty}$ sense. Indeed, it is a generalized cross-cap, that is, $\mathscr{A}$--equivalent to the map of \cref{eq:localWU}.
\end{enumerate}
\end{proposition}

\begin{remark}\label{re:carriesthechange}
    Notice that, by \cref{it1:propertiesofshath,it3:propertiesofshath}, the boundary of the image of $\shat{h}_{\theta}$ decomposes as $(B^2 \vee B^2) \cup_{\partial} (\Sbb^1 \times [-1,1])$, hence $\shat{h}_\theta$ carries the change from  the fibers of $K_f \cap T_{f\{i\}}$ to the fibers of $F_f \cap T_{f\{i\}}$, cf. \cref{lem:twORuntw} and \cref{eq:TC2}. Later we use this fact, as well as the properties listed above.
\end{remark}

\color{black}

Finally, we define $\shat{h}: X_{\textnormal{ext}} \to \shat{B}_{\textnormal{ext}}$. For each untwisted component we take
\begin{equation*}
\begin{aligned}
    	\shat{h}|:\Sbb^1\times B^2\times \left[-1,1\right]&\longrightarrow \Sbb^1\times B^4 \times \left[-1,1\right]\\
    \shat{h}(\theta;\: z;\:t)&= \big(\theta; \: \shat{h}_\theta(z;\: t)\big).
\end{aligned}
\end{equation*}
For twisted components we take
\begin{equation*}
    \shat{h}|:Z= \frac{\left[-\pi,\pi\right]\times B^2\times\left[-1,1\right]}{(-\pi,z,t)\sim(\pi,\overline{z},-t)} \to 
    \frac{[-\pi, \pi] \times B^4 \times \left[-1,1\right]}
    {(-\pi;\: \alpha, \beta;\: s) \sim (\pi;\: \beta, \alpha;\: s)}
\end{equation*}
defined, locally, in the same way as we did for untwisted components. While it is easy to see that it is a well-defined $\mathcal{C}^{\infty}$ map for untwisted components, it is not obvious for twisted components. Recall that the target space of $ \shat{h}|$ is always diffeomorphic to $\Sbb^1 \times B^4 \times [-1, 1]$.

\begin{proposition} For twisted components,  $ \shat{h}|$ is a well-defined $\mathcal{C}^{\infty}$ map.
\end{proposition}
\begin{proof} We have the following situation:
\begin{equation*}
 \
\begin{tikzcd}
	{\left[-\pi,\pi\right]\times B^2\times\left[-1,1\right]} & {\frac{\left[-\pi,\pi\right]\times B^2\times\left[-1,1\right]}{(-\pi;\;z;\;t)\sim(\pi;\;\overline{z};\;-t)}} \\
	{\left[-\pi,\pi\right]\times B^4 \times \left[-1,1\right]} & { \frac{[-\pi, \pi] \times B^4 \times \left[-1,1\right]}
    {(-\pi;\; \alpha, \beta;\; s) \sim (\pi;\; \beta, \alpha;\; s)}}
	\arrow["q_1", from=1-1, to=1-2]
	\arrow["{\id\times\shat{h}_\theta}", from=1-1, to=2-1]
	\arrow["q_2", from=2-1, to=2-2]
	\arrow["{\shat{h}}", from=1-2, to=2-2]
\end{tikzcd},\end{equation*}
where $q_1$ and $q_2$ are the quotient maps defining the bundles over $\Sbb^1$.
     We have to check that the map commutes with the monodromy of the source and the target. That is, 
    \begin{equation*}\shat{h} (-\pi, z, t) \textnormal{ and } \shat{h} (\pi, \overline{z}, -t)
    \end{equation*}
    glue together by the monodromy of the target space, i.e. by swapping the the complex coordinates of $B^4$. This holds, since
    \begin{equation*}
   \begin{aligned}
    \shat{h} (-\pi, z, t)=&
     \left(-\pi;\:
 \frac{z(1+t)}{\sqrt{2+2t^2}}, \ 
 \frac{\overline{z}(1-t)}{\sqrt{2+2t^2}}
 ;\:-t^2 \right)\\
     \shat{h} (\pi, \overline{z}, -t)=&
     \left(\pi;\:
 \frac{\overline{z}(1-t)}{\sqrt{2+2t^2}}, \ 
 \frac{z(1+t)}{\sqrt{2+2t^2}}
 ;\:-t^2 \right).
\end{aligned}
    \end{equation*}
\end{proof}

\subsubsection{The gluing of the interior and the exterior maps}

\begin{lemma}\label{lem:h01}
There is a $\mathcal{C}^{\infty}$ map $h:X\to\shat{B}$  such that 
\begin{enumerate}[label={(\roman*)}]
	\item $h$ restricted to the interior part $B^4\subset X$ coincides with $f_s^\RR$,\label{it:1h1}
   \item $h$ restricted to $ X_{\textnormal{ext}}$ coincides with  $\shat{h}$. 
\end{enumerate}
\end{lemma}
\begin{proof}
The situation is given by the following diagram, where $Z_{f\{i\}}$ is given in \cref{eq:Ziinsecondform}.
\[\begin{tikzcd}[column sep=12pt]
	&& {N_i} && {Z_{f\{i\}}} \\
	X && {B^4} && {X_{\textnormal{ext}}} \\
	{\shat{B}} && {B^6} && {\shat{B}_{\textnormal{ext}}} \\
	&& {T_{f\{i\}}} && {\mathbb{S}^1\times B^4\times \{-1\}}
	\arrow["h"', dashed, from=2-1, to=3-1]
	\arrow["{f_s^\RR}"', from=2-3, to=3-3]
	\arrow["{\shat{h}}", from=2-5, to=3-5]
	\arrow[hook, from=1-5, to=2-5]
	\arrow[hook, from=4-5, to=3-5]
	\arrow["\Phi", from=1-3, to=1-5]
	\arrow["{\phi^{f\{i\}}}"', from=4-3, to=4-5]
	\arrow[hook, from=4-3, to=3-3]
	\arrow[hook, from=1-3, to=2-3]
	\arrow["{\bigcup_{\mbox{$\Phi$}}}"{description}, draw=none, from=2-3, to=2-5]
	\arrow["{\bigcup_{\mbox{$\varphi$}}}"'{description}, draw=none, from=3-3, to=3-5]
	\arrow["="{description}, draw=none, from=2-1, to=2-3]
	\arrow["{=}"{description}, draw=none, from=3-1, to=3-3]
\end{tikzcd}\]

To obtain ${h}$, we need to check that the gluings by $\Phi$ and $\varphi$ are compatible with $f_s^\RR$ and $\shat{h}$. 

By \cref{sss:targetspace,eq:phiT1},  the gluing $\varphi$ defining $\shat{B}$ identifies $K_f\cap T_{f\{i\}}$ with a locally trivial fiber bundle with fiber $\{\alpha \beta=0\}$ in $\Sbb^1 \times B^4 \times \{-1\}$. By \cref{it3:propertiesofshath} of \cref{pr:propertiesofshath}, this agrees set-wise with the image of $\shat{h}$. However, this leaves some freedom in the choice of the diffeomorphisms $\phi^{f\{i\}}$ that make $T_{f\{i\}} \cong \Sbb^1 \times B^4$. More precisely, they could be composed by a diffeomorphism fixing (set-wise) $\{\alpha=0\}$ and $\{\beta=0\}$. In contrast, $X$ is already determined (by $\Phi$), so these $\phi^{f\{i\}}$ have to yield something compatible with an image of $X$ (by both $f_s^\RR$ and $\shat{h}$).

Now, we determine this freedom of $\phi^{f\{i\}}$ in terms of the fixed diffeomorphism $\Phi$ and the maps $f_s^\RR$ and $\shat{h}$. This is done by the diagram chasing of the points:
\begin{equation}\label{eq:diagramchasing}
\begin{tikzcd}
	{(f_s^\RR)^{-1}(p_1)} & {(z,+1)} & {(z,-1)} & {(f_s^\RR)^{-1}(p_2)} \\
	{p_1} & {(\alpha,0)} & {(0,\beta)} & {p_2}
	\arrow["\Phi", maps to, from=1-1, to=1-2]
	\arrow[dashed, leftrightarrow, from=1-2, to=1-3]
	\arrow["\Phi"', maps to, from=1-4, to=1-3]
	\arrow["{\shat{h}}", maps to,from=1-2, to=2-2]
	\arrow["{(f_s^\RR)^{-1}}", maps to, from=2-1, to=1-1]
	\arrow["{\phi_{\alpha}}"', maps to, from=2-1, to=2-2]
	\arrow["{(f_s^\RR)^{-1}}"', from=2-4, to=1-4]
	\arrow["{\shat{h}}", maps to, from=1-3, to=2-3]
	\arrow["{\phi_{\beta}}", maps to, from=2-4, to=2-3]
\end{tikzcd}
\end{equation}
For a fixed $\theta \in f(\delta_i)$, take a point $p_1 \in T_{f\{i\}} \cap K_f \subset \Sbb^5$ on the local branch to be identified with $\{\beta=0\}$ by $\phi^{f\{i\}}$. Then, $p_1$ has a unique preimage by $f_s^{\R}$ in $N_i \subset \Sbb^3$, which glues by $\Phi^{f\{i\}}$ to a point $(z, +1) \in B^2 \times [-1,1]$, as a part of either $\Sbb^1 \times B^2 \times [-1, 1]$ or $Z$ ($f\{\delta_i\}$ untwisted or twisted, respectively). Then, $\phi(p_1)$ has to be defined as $\shat{h}_{\theta}(z,+1)$. 

However, $p_1$ has a pair $p_2$ on the other local branch that is identified with $\{\alpha=0\}$, and $\phi(p_2)$ has to be defined as $\shat{h}_{\theta}(z, -1)$. Hence, $\phi^{f\{i\}}$ is fixed simultaneously on the two branches corresponding to $\{\beta=0\}$ and $\{\alpha=0\}$, i.e., the restrictions determine each other  (denoted by $\phi_{\alpha}$ and $\phi_{\beta}$, respectively).
This is only possible if the left and right squares of the diagrams are compatible by the symmetry between the points $(z,1)$ and $(z,-1)$. This can only be obstructed by orientation. 

On the one hand, $\phi^{f\{i\}}$ has to be an orientation-preserving diffeomorphism from $T_{f\{i\}}$ to $\Sbb^1 \times B^4 $, to obtain  an oriented manifold $\shat{B}$ by the gluing $T_{f\{i\}} \to \Sbb^1 \times B^4 \times \{-1\} $. Here, $\Sbb^1 \times B^4 $ is considered with its orientation given by $(\theta, \alpha, \beta)$, but $\Sbb^1 \times B^4 \times \{-1\}$ inherits the opposite orientation as a part of the boundary $ \partial (\Sbb^1 \times B^4 \times [-1, 1])$. To obtain an orientation-preserving $\phi^{f\{i\}}$, the diffeomororphisms $\phi_{\alpha}$ and $\phi_{\beta}$ between complex discs have to preserve or reverse the orientation simultaneously.

On the other hand, we analyse the maps $f_s^{\R}$, $\shat{h}$ and $\Phi$, concluding that the above construction gives both $\phi_{\alpha}$ and $\phi_{\beta}$ orientation reversing (cf. \cref{re:oribound}). The map $f_s^{\R}$ (equal to the holomorphic map $f$ near the boundary) preserves the orientation between the disc-fibers of $N_i$ (and $N_{i'}$) and the fibers $B^2 \vee B^2$ of $K_f \cap T_{f\{i\}}$. The  map $\shat{h}_{\theta}$ preserves the orientation for $B^2 \times \{1\} \to \{\beta =0\}$, but it reverses the orientation for $B^2 \times \{-1\} \to \{\alpha =0\}$, see \cref{it3:propertiesofshath ori} of \cref{pr:propertiesofshath}. The gluing $\Phi^{f\{i\}}$ (for an untwisted component, twisted components are similar)  preserves the orientation between $B^2 \times \{-1\}$ and the disc-fibers of $N_i$, and it reverses the orientation between $B^2 \times \{1\}$ and the disc-fibers of $N_{i'}$, because $X$ becomes an oriented manifold only in this way. Indeed, both discs $B^2 \times \{-1\}$ and $B^2 \times \{1\}$ are endowed with their complex orientations, which agree with the  orientations inherited from $\partial (B^2 \times [-1,1])$ for $B^2 \times \{1\}$, and they are the opposite of that for $B^2 \times \{-1\}$.

Putting these together, both $\phi_{\alpha}$ and $\phi_{\beta}$ reverse the orientation by the above construction (see the diagram \cref{eq:diagramchasing}). For $\phi_{\alpha}$, the orientation is reversed by $\Phi$ (and it is preserved by $f_s^{\R}$ and $\shat{h}$) and, for $\phi_{\alpha}$, the orientation is reversed by $\shat{h}$ (and it is preserved by $\Phi$ and $f_s^{\R}$).
This gives the correct orientation, makes $\phi$ well-defined and compatible with the diagram. 
\end{proof}

\subsection{Properties of the map}\label{ssec:properties}
In this section we show that the map $h$ from \cref{lem:h01} provides a slice singular manifold for $\iota|: \partial F_f \hookrightarrow \Sbb^5$ and we compute the terms of the ESzTS formula \cref{th:saekiszucs}.

Without further analysis, it is easy to show the following.
\begin{proposition}
By construction, $h$  as defined in \cref{lem:h01} is stable in the $\Ccal^\infty$ sense.
\end{proposition}

\subsubsection{The map $h$ on the boundary}\label{sss:themapontheboundary}
Recall  \cref{th:milnperpart1} for the construction of $\partial F$. By  \cref{prop:xcongf},  $\partial X\cong \partial F_f$. Moreover, the stable map $h: X\to\shat{B}\cong B^6$, defined by \cref{lem:h01}, induces an embedding  $h|_{\partial X}: \partial X\to \partial \shat{B}\cong \Sbb^5$, cf. \cref{it1:propertiesofshath} of \cref{pr:propertiesofshath}.

However, this is not enough to guarantee that $h$ is a slice singular manifold for $\partial F$ and, later, to use the ESzTS formula \cref{th:saekiszucs} in the way we want. We need that $h|_{\partial X}$ is regular homotopic to the inclusion $\iota|:\partial F\hookrightarrow\Sbb^5$ in $\Sbb^5$. Indeed, we show that this construction of $h$ makes them isotopic.

There is a slight imprecision to solve. The target spaces of $h|_{\partial X}$ and $\iota|$ are not the same, they are only diffeomorphic. In order to compare both embeddings we have to fix a diffeomorphism $\Xi_\partial: \partial \shat{B} \to \Sbb^5=\partial B^6$ as follows. 

First, observe that there is a diffeomorphism
\begin{equation*}
\Xi_{\textnormal{cyl}}:\left(\partial B^4 \times [-1,1]\right) \cup_{\partial} \left(B^4 \times \{1\} \right) \rightarrow B^4 \times \{-1\}
\end{equation*}
so that $\Xi_{\textnormal{cyl}}$ takes $B^4 \times \{1\}$ to a small ball $\bm{B}\subset B^4\times \{-1\}$ around zero and $\partial B^4 \times [-1,1]$ to $\text{cl}\big(B^4 \times \{-1\} \setminus\bm{B}\big)$ extending $\bm{B}$ in a radial way (see \cref{fig:cylinder}).

Then, since 
$$\shat{B}_{\textnormal{ext}}=\bigsqcup_{f\{i\}} \Sbb^1\times B^4 \times \left[-1,1\right]$$
is glued to $B^6$  along $\Sbb^1\times B^4 \times \left\{-1\right\}$ (see \cref{sss:targetspace}), we take an extension $\Xi_\partial$ of $\Xi_{\textnormal{cyl}}$ to $\partial \shat{B}$ in the obvious way along the $\Sbb^1$-families
$$ \bigsqcup_{f\{i\}} \Sbb^1\times \left(\left(\partial B^4 \times [-1,1]\right)  \cup_{\partial} \big(B^4 \times \{1\}\big)\right).$$
It is also easy to see that there is another diffeomorphism $\Xi:\shat{B} \to B^6$ extending $\Xi_\partial$, since we are simply pushing down the exterior part $\shat{B}_{\textnormal{ext}}$ to $B^6$.

\begin{theorem}\label{thm:reghom of h}
The embeddings
$$\begin{aligned}
\Xi_\partial \circ h|_{\partial X}:\partial X & \hookrightarrow \partial \shat{B} \to \Sbb^5\\
\iota|:\partial F_f &\hookrightarrow \Sbb^5
\end{aligned} $$
as defined above are isotopic in $\Sbb^5$. Consequently, $\Xi\circ h: X \to B^6$ is a slice singular manifold of $\iota|: \partial F \hookrightarrow \Sbb^5$.
\end{theorem}

First, let us study a local prototype of the theorem. 
For simplicity, we denote  $B^4 \times \{-1\}$ by $B^4 $. Observe that the Hopf link 
\begin{equation}\label{eq:hopflink}
    \{\alpha \beta=0\} \cap \Sbb^3 \subset \partial B^4
\end{equation}
has the Milnor fiber ($\delta > 0$)
\begin{equation}\label{eq:slice1}
    \{\alpha \beta=\delta\} \subset B^4 
\end{equation}
as a slice surface after a canonical isotopy on the boundary. By \cref{it1:propertiesofshath,it3:propertiesofshath} of \cref{pr:propertiesofshath}, another slice surface is given by the embedding
\begin{equation}\label{eq:slice2prevform}
   \begin{tikzcd}
	{\Sbb^1 \times [-1, 1]} & {\Sbb^3 \times [-1, 1]} & {B^4 \setminus \bm{B}}
	\arrow["{\shat{h}_{\theta}|}", hook, from=1-1, to=1-2]
	\arrow["{\Xi_{\text{cyl}}}", hook, from=1-2, to=1-3]
\end{tikzcd}.
\end{equation}
\begin{lemma}\label{le:slicesurf}
    The slice surfaces of the Hopf-link \cref{eq:hopflink} given by the Milnor fiber and $\Xi_{\textnormal{cyl}} \circ    \shat{h}_{\theta}|$ (\cref{eq:slice1,eq:slice2prevform}) are isotopic in $B^4$ by an isotopy that fixes the Hopf link point-wise.
\end{lemma}
\begin{proof}[Proof of \cref{le:slicesurf}]
    Consider the radial homotopy $p: \left(B^4 \setminus \{0\}\right) \times\left[-1,1\right] \rightarrow B^4 $,
    \begin{equation*}
         p(sv, sw;\: t)= p_t(sv,sw)=\big((s+(1-s)t)v,(s+(1-s)t)w\big),
    \end{equation*}
where $(v, w) \in \Sbb^3$, $0<s \leq 1$.

First, we show that
\begin{enumerate}
    \item\label{it0slicesurf} the restriction of the homotopy $p_t$ to the slice surfaces is an isotopy;
    \item\label{it1slicesurf} the restriction of $p_1$ to each of the two slice surfaces is an embedding into $\Sbb^3$ (hence, a Seifert surface of the Hopf link); and
    \item\label{it2slicesurf} their images by $p_1$ define the same Seifert surface.
\end{enumerate}

\cref{it0slicesurf,it1slicesurf} are standard arguments for both slice surfaces (recall the definition of $\shat{h}_\theta$ \cref{eq:w,eq:defhtheta}). The image of the first slice surface is
\begin{equation*}
    p_1\big(\{(\alpha, \beta) \in B^4 \ | \ \alpha \beta=\delta\}\big)=\big\{(v,w) \in \Sbb^3 \ | \ vw \in \R^+ \big\}.
\end{equation*}
Indeed, for $\alpha=sv$ and $\beta=sw$, we get $vw=\delta/s^2$, so $vw \in \R^+$ and $vw \geq \delta$. Taking the limit $\delta \to 0$, the boundary of this surface is exactly the Hopf link $\{vw=0\}$. Notice also that the maximal value of $|vw|$ on $\Sbb^3$ is $1/2$.

For the second slice surface, observe that 
\begin{equation*}
    p_1 \circ \Xi_{\textnormal{cyl}} \circ    \shat{h}_{\theta}| = \pi_1 \circ    \shat{h}_{\theta}|,
\end{equation*}
where $\pi_1: \Sbb^3 \times [-1, 1] \to \Sbb^3$ is the projection onto the $\Sbb^3$ component. Hence, we have to consider the first two coordinates in \cref{eq:w}: $(v, w) \in \Sbb^3$ is in the image of the projection if, and only if, there is a $z \in \Sbb^1$ and $t \in [-1, 1]$ such that
\begin{equation}\label{eq:vw}
   (v,w)= \left(\frac{z(1+t)}{\sqrt{2+2t^2}},\frac{\overline{z}(1-t)}{\sqrt{2+2t^2}}\right).
\end{equation}
Therefore, 
\begin{equation}\label{eq:vwt}
    vw=\frac{1-t^2}{2+2t^2},
\end{equation}
which implies that $vw \in \R^+$ and $0 \leq vw \leq 1/2$. The converse is also true: for any point $(v,w) \in \Sbb^3$ with $0 \leq vw \leq 1/2$ there is a unique $z \in \Sbb^1$ and $t \in [-1, 1]$ satisfying \cref{eq:vw}. Indeed, 
the right side of \cref{eq:vwt} expressing $vw$ determines $t$ up to sign, and by \cref{eq:vw} the sign of $t$ has to agree with the sign of $|v|-|w|$.

Since both Seifert surfaces obtained by the projections  to $\Sbb^3$ of the two slice surfaces coincide, the slice surfaces are isotopic in $B^4$ by the radial homotopy $p_t$.
\end{proof}

\begin{proof}[Proof of \cref{thm:reghom of h}]
    The map $\Xi_\partial \circ h|_{\partial X}$ and its parts is summarized in the following diagram. 
\[\begin{tikzcd}
	{\partial B^4\setminus\partial X_\text{ext}} & {\partial B^6\setminus\partial\shat{B}_\text{ext}} & {\partial B^6} \\
	{\partial X} & {\partial\shat{B}} & {\partial B^6} \\
	{\partial X_\text{ext}\setminus\partial B^4} & {\partial \shat{B}_\text{ext}\setminus\partial B^6} & {\partial B^6} \\
	{\partial B^2 \times[-1,1]} & {\partial B^4\times[-1,1]} & {B^4\setminus\bm{B}}
	\arrow[hook, from=4-1, to=3-1]
	\arrow[hook, from=1-1, to=2-1]
	\arrow[hook, from=3-1, to=2-1]
	\arrow[equal, from=3-3, to=2-3]
	\arrow[equal, from=1-3, to=2-3]
	\arrow["{\shat{h}_\theta|}", from=4-1, to=4-2]
	\arrow["{f_s^\RR|}", from=1-1, to=1-2]
	\arrow["{h|}", from=2-1, to=2-2]
	\arrow["{\shat{h}|}", from=3-1, to=3-2]
	\arrow["{\Xi_\partial}", from=2-2, to=2-3]
	\arrow["{\Xi_\partial|}", from=3-2, to=3-3]
	\arrow[hook, from=1-2, to=1-3]
	\arrow["{\Xi_\text{cyl}|}", from=4-2, to=4-3]
	\arrow[hook, from=4-3, to=3-3]
	\arrow[hook, from=4-2, to=3-2]
	\arrow[hook, from=1-2, to=2-2]
	\arrow[hook, from=3-2, to=2-2]
\end{tikzcd}\]

The second row is the map we want to study, $\Xi_\partial\circ h|_{\partial X}$.

The first row is the restriction of $\Xi_\partial\circ h|_{\partial X}$ to $\Sbb^3 \setminus \bm{N}$,
 which coincides with $f_s^\RR$ (and also $f$). Its image is isotopic to $\partial F_f \setminus \bm{T} \subset \Sbb^5 \setminus \bm{T}$. Indeed, this corresponds to the non-singular part of $K_f$, i.e., $K_f\setminus \bm{T}$ (recall \cref{sss:surgerypartialf part1}).


The third row is the restriction of $\Xi_\partial\circ h|_{\partial X}$ to the boundary part of $\partial X$ corresponding to the exterior part (that is, $\partial X_\text{ext}\setminus\partial B^4$). This map is a collection of  $\Sbb^1$-families of embeddings $\Xi_{\text{cyl}}\circ \shat{h}_\theta|$ as in \cref{eq:slice2prevform}, i.e., the last row.
Furthermore, $\bm{T}$ is mapped by $\varphi$ to the $\Sbb^1$-families of $B^4 \times \{-1\}$ (see \cref{sss:targetspace} and, specifically, \cref{eq:phiT1}).

Finally, the inclusion of $F_f\cap \bm{T}$ by $\varphi$ is isotopic to an $\Sbb^1$-family of the Milnor fibers $\{\alpha\beta=\delta\}$, cf. \cref{ss:surgerypartialf} and in particular \cref{eq:TC2}. By \cref{le:slicesurf}, this is isotopic to the restriction of $\Xi_\partial\circ h|_{\partial X}$ to the boundary part of $\partial X$ corresponding to the exterior part, by an isotopy that fixes the $\Sbb^1$-families of Hopf links (i.e., the intersection between the exterior and the interior parts).

The result follows, since we have an isotopy on each part and the intersection is fixed.
\end{proof}

 \begin{remark}[On the orientation]\label{re:oribound}
Although we proved \cref{le:slicesurf}  by a direct computation (and \cref{thm:reghom of h} through \cref{le:slicesurf}), the orientation plays a very important role in our construction. It is hidden in the form of the generalized real cross-cap we use: we changed the orientation on one disc, see \cref{eq:w} and \cref{it3:propertiesofshath ori} of \cref{pr:propertiesofshath}. 
See also \cref{sss:targetspace} and the proof of \cref{lem:h01}.

 \end{remark}

\subsubsection{The linking numbers $\ell$ and $L$}\label{sss:linkinglL}

We characterize now the linking numbers $\ell$ and $L$ in the ESzTS formula \cref{th:saekiszucs} for our slice singular manifold $h: X \to \shat{B}$ given in \cref{ss:constructionh,ss:4mfld}.

$L(h|_{\partial{X}})$ is very easy since ${h}|_{\partial X}: \partial X \hookrightarrow \partial \shat{B}$ is an embedding, hence $L(h|_{\partial X})=0$ by definition of $L$ (see \cref{def:L}).

For $\ell$, we prove the following.

\begin{theorem}\label{thm:l=c}
For ${h}:X\to\shat{B}$ given in \cref{ss:constructionh,ss:4mfld},
    $$\ell({h})=C(f).$$
   \end{theorem}

Recall $\ell(h)$ as defined in \cref{def:ell}. Let $\mathcal{W} \subset \shat{B}$ denote the set of non-immersive points (singular points) of the generic map ${h}$ in the target. They are generalised Whitney umbrella points (cf. \cref{eq:localWU}) and they are contained in the boundary of the double values of ${h}$. Let $\mathcal{W}' \subset \shat{B}$ denote the copy of $\mathcal{W}$ shifted  slightly along the outward normal field of $\mathcal{W}$, as the boundary of the double values. Then,
\begin{equation*}
    \ell(h)\coloneqq \mbox{lk}_{({\scriptsize \shat{B}},{\scriptsize\partial\shat{B}})}\big(\mathcal{W}', h(X)\big).
\end{equation*}

The map ${h}$ has two types of generalised Whitney umbrella points:
\begin{enumerate}
    \item The singular points of $h|_{B^4}$, i.e. the real $\mathcal{C}^{\infty}$ deformation $f_s^\RR$ of $f$ inside $B^6$. We call them \textit{interior singular points} of ${h}$, denoted by $\mathcal{W}_{\mathrm{int}}$. Its components correspond to the $\mathcal{C}^{\infty}$-stabilization of the complex Whitney umbrella points of $f_s$, given in \cref{ss:assoc}.
    \item The \textit{exterior singular points}, i.e., the singular points of the exterior map  $\shat{h}: X_{\textnormal{ext}} \to \shat{B}_{\textnormal{ext}}$ in the target. In the source, the locus of singular points is the union of the middle circles of the gluing pieces, cf. \cref{eq:whuntw} and \cref{eq:whtw}. Its image is
    \begin{equation*}
\mathcal{W}_{\mathrm{ext}}\coloneqq \bigsqcup_{f\{i\}} \Sbb^1 \times \{0\} \times \{0\}\subset \bigsqcup_{f\{i\}} \Sbb^1\times B^4\times\left[-1,1\right],
    \end{equation*}
cf. \cref{it56:propertiesofshath}  of \cref{pr:propertiesofshath}.
\end{enumerate}

According to the decomposition $\mathcal{W}= \mathcal{W}_{\mathrm{int}}\sqcup\mathcal{W}_{\mathrm{ext}}$, with the obvious notation,
\begin{equation*}
    \ell({h})=\ell_{\mathrm{int}}({h}) + \ell_{\mathrm{ext}}({h}).
\end{equation*}

\begin{proposition}
For ${h}:X\to\shat{B}$ given in \cref{ss:constructionh,ss:4mfld},
$$\ell_{\mathrm{int}}({h})=C(f).$$
\end{proposition}
\begin{proof}
    It is proved in \cite[Theorem 9.1.3]{Nemethi2015} that for the $\mathcal{C}^{\infty}$-stabilization of the complex Whitney umbrella, given in \cref{ss:assoc}, $\ell=+1$. Cf. \cref{pr:summ}.
\end{proof}

\begin{theorem}\label{th:lext}
For ${h}:X\to\shat{B}$ given in \cref{ss:constructionh,ss:4mfld},
$$\ell_{\mathrm{ext}}({h})=0.$$
\end{theorem}
\begin{proof}
     By definition, 
     \begin{equation*}
\ell_{\mathrm{ext}}({h})=\text{lk}_{({\scriptsize \shat{B}},{\scriptsize\partial\shat{B}})}\big(\mathcal{W}_{\text{ext}}',{h}(X)\big),
     \end{equation*}
for a shifted copy $\mathcal{W}_{\text{ext}}'$ of $\mathcal{W}_{\text{ext}}$.

$\mathcal{W}_{\mathrm{ext}}$ decomposes as $\bigsqcup_{f\{i\}} \mathcal{W}_{f\{i\}}$ according to the components of $f(\bm{\delta})$. We prove that the linking number $\ell_{f\{i\}}({h})$ corresponding to every part $f(\delta_i)$ is 0, where 
\begin{equation}
\ell_{f\{i\}}({h})\coloneqq\text{lk}_{({\scriptsize \shat{B}},{\scriptsize\partial\shat{B}})}\big(\mathcal{W}_{f\{i\}}',{h}(X)\big). 
     \end{equation}
A shifted copy $\mathcal{W}'_{f\{i\}}$ can be given explicitly as
\begin{equation}\label{eq:pushing}
\mathcal{W}'_{f\{i\}}=\Sbb^1 \times \{0\} \times \{\epsilon\}\subset \Sbb^1\times B^4\times\left[-1,1\right],
    \end{equation}    
    for a small positive number $\epsilon$ (see \cref{fig:cylinder}). Indeed,  the set of double points of $\shat{h}$ is $\Sbb^1 \times \{0\} \times [-1, 0]$ (cf. \cref{pr:propertiesofshath}), so $\mathcal{W}'_{f\{i\}}$ is a shifted copy of $\mathcal{W}_{f\{i\}}$ along its outward normal field, as a boundary of the double points.
     
Let us define $\mathcal{W}''_{f\{i\}}$ by pushing out $\mathcal{W}'_{f\{i\}}$ further along the normal field, until it reaches the boundary (see \cref{fig:cylinder}), that is,
\begin{equation*}
\mathcal{W}''_{f\{i\}}=\Sbb^1 \times \{0\} \times \{1\}\subset \left(\Sbb^1\times B^4\times\left[-1,1\right] \right) \cap \partial\shat{B}.
    \end{equation*}    

On the one hand, we show that 
    \begin{equation*}
\text{lk}_{({\scriptsize \shat{B}},{\scriptsize\partial\shat{B}})}\big(\mathcal{W}_{f\{i\}}',{h}(X)\big)=
\pm \text{lk}_{{\scriptsize\partial\shat{B}}}\big(\mathcal{W}_{f\{i\}}'',{h}(\partial X)\big).
     \end{equation*}
    Indeed, by \cref{re:linkingdef},
\begin{equation*}
\text{lk}_{({\scriptsize \shat{B}},{\scriptsize\partial\shat{B}})}\big(\mathcal{W}_{f\{i\}}',{h}(X)\big)=
\text{int}_{({\scriptsize \shat{B}},{\scriptsize\partial\shat{B}})}\big(\Lcal',{h}(X)\big),
     \end{equation*}
     where $\Lcal'$ is a homological membrane (a surface) contained in $\shat{B}$ so that $\partial \Lcal'=\mathcal{W}_{f\{i\}}'$. Furthermore, $\Lcal'$ can be extended with the collar $\Sbb^1\times \{0\}\times\left[\epsilon,1\right]$ to obtain the membrane $\Lcal''$, that has $\partial \Lcal''=\mathcal{W}_{f\{i\}}''$. Then, 
     \begin{align*}
        \text{int}_{({\scriptsize \shat{B}},{\scriptsize\partial\shat{B}})}\big(\Lcal',{h}(X)\big)&=\text{int}_{({\scriptsize \shat{B}},{\scriptsize\partial\shat{B}})}\big(\Lcal'',{h}(X)\big)\\
     &= \pm \text{lk}_{{\scriptsize\partial\shat{B}}}\big(\mathcal{W}_{f\{i\}}'',{h}(\partial X)\big).
     \end{align*}

     Note that we used two different definitions of the linking number, see \cref{re:linkingdef}.

     On the other hand, we show that
     \begin{equation*}
        \text{lk}_{\scriptsize\partial\shat{B}}\big(\mathcal{W}_{f\{i\}}'',{h}(\partial X)\big)=0.
     \end{equation*}
To do so, we use the diffeomorphism $\Xi_{\partial}: \partial\shat{B} \to \Sbb^5$  to copy the whole configuration from $\partial\shat{B}$ to $\Sbb^5$. The composition of $\Xi_\partial$ with the isotopy indicated by \cref{thm:reghom of h} takes $h(\partial X)$ to $\partial F_f $ and $\Wcal_{f\{i\}}''$ to  $f(\delta_i)$.
Hence,
\begin{equation}\label{eq:lmilnordouble}
\text{lk}_{\scriptsize\partial\shat{B}}\big(\mathcal{W}_{f\{i\}}'',h(\partial X)\big)=
\text{lk}_{\Sbb^5}\big(f(\delta_i),\partial F_f\big)=0,
\end{equation}
 since $f(\delta_i)=\partial f(D_i)$, and $f(D_i) \cap F_f=\varnothing$, because they are on different level sets of $g$.
\end{proof}

\begin{remark}\label{re:linkingdef} We used two different definitions of the linking number, we recall them here. Let $A, B  \subset \Sbb^n$ be two homological cycles with $\Z$ coefficients (e.g., images of oriented closed manifolds by continuous maps) of dimensions $0<a$ and $0<b$ with $a+b=n-1$. Assume that $A \cap B= \varnothing$. Let $\Lcal^A, \Lcal^B$ be `slice manifolds' of $A$ and $B$, that is, two chains  in $B^{n+1}$ of dimensions $a+1$ and $b+1$ such that $\partial \Lcal^A=A$ and $\partial \Lcal^B=B$. Let $\Scal^A $ be a `Seifert manifold' of $A$, that is, a chain in $\Sbb^n$ of dimension $a+1$ with $\partial \Scal^A=A$. Then, by both definitions,
\begin{equation}\label{eq:linkingnumber}
    \mbox{lk}_{\Sbb^n} (A, B)=  \mbox{int}_{B^{n+1}}(\Lcal^A, \Lcal^B)= \pm
    \mbox{int}_{\Sbb^n}(\Scal^A, {B}),
\end{equation}
where $\mbox{int}$ denotes the intersection number, that is, the algebraic number of the transverse intersection points (the intersection of the chains is assumed to be transverse). The correct sign might depend on the convention.

To see the equality of the two intersection numbers up to sign, consider $B^{n+1}$ embedded in $ \Sbb^{n+1}$ as the northern hemisphere. Let $ \hat{B}$ the mirror image of $\Lcal^B$, reflected by the equator hyperplane. Then $\Lcal^B + \hat{B}$ and $\Lcal^A - \Scal^A$ are closed chains (cycles) in $\Sbb^{n+1}$ of dimension less than $n+1$, hence
\begin{equation*}
    \mbox{int}_{\Sbb^{n+1}}( \Lcal^A - \Scal^A, \Lcal^B + \hat{B})=0.
\end{equation*}
For the terms appearing on the left-hand side we have 
\begin{equation*}
\begin{aligned}
\mbox{int}_{\Sbb^{n+1}}(\Lcal^A , \Lcal^B )=&\mbox{int}_{B^{n+1}}(\Lcal^A , \Lcal^B ),\\
\mbox{int}_{\Sbb^{n+1}}(\Lcal^A , \hat{B} )=&0,\\
\mbox{int}_{\Sbb^{n+1}}(  - \Scal^A, \Lcal^B + \hat{B})=&
    \pm \mbox{int}_{\Sbb^{n}}(   \Scal^A, B),
\end{aligned}
\end{equation*}
proving the equality of the intersection numbers up to sign in \cref{eq:linkingnumber}. The sign might change when we step back from $\Sbb^{n+1}$ to $\Sbb^n$ in the last equation, depending on the convention. However, in the proof of \cref{th:lext} it is zero, so this is irrelevant.
\end{remark}

\begin{remark}
As the trivial local description of the pushing out \cref{eq:pushing} suggests, \cref{th:lext} has a global nature, indeed, it is based on the global relation of the boundary of the Milnor fiber and the double point curve of $f$, cf. \cref{eq:lmilnordouble}. This property cannot seem locally, it is hidden in the global picture given by  the vertical indices.

Comparing it with \cite[Proposition 1.1 and Section 4]{Pinter2023} is interesting, in which `nearby embedded 3-manifolds' of the associated immersion $f|_{\Sbb^3}$ are defined, and $\partial F_f$ is characterised as the one which has zero linking number with the double point curve. The whole package of the present paper can be generalized to the nearby 3-manifolds, a slice singular manifold can be defined for it in the same way we did  for $\partial F_f$, and by \cref{th:saekiszucs} a regular homotopy invariant can be associated to them from the properties of the slice singular manifold. In this case, the exterior linking number will be not zero in general.
\end{remark}

\subsubsection{The algebraic number of triple values}

From \cref{sec:ssmoi} , it is very easy to see the following.
\begin{proposition}\label{prop:t=T}
For the constructions of $h:X\to\shat{B}$ given in \cref{ss:constructionh,ss:4mfld},
   $$t(h)= T(f).$$
\end{proposition}
\begin{proof}
This follows from the construction of $h$: the only triple values come from $f_s^\RR$ and we have $t(f_s^{\R})=T(f)$, see \cref{pr:summ} (cf. \cite{Nemethi2015}).  Indeed, no new triple values are created in the real stabilization to go from $f_s$ to $f_s^\RR$ nor in the extension to go from $f_s^\RR$ to $h$. Since $f_s$ is complex, every triple value contributes with $+1$ to the total sum of $t(h)$.
\end{proof}

\section{Main theorem} \label{sec:main}

We can now prove our main theorem. Recall that first we consider a finitely generated map germ $f:(\CC^2,0)\to(\CC^3,0)$, then we consider the surface $\big(\im(f),0\big)$ with reduced structure and, finally, we take its Milnor fiber $F_f$ and compute its signature $\sigma(F_f)$. We do it by means of the ESzTS formula \cref{th:saekiszucs}; using the construction of $X$ in \cref{ss:4mfld}; and, for the proof, a slice singular manifold $h:X\to\shat{B}\cong B^6$ for $\partial F_f$ given in \cref{ss:constructionh}. 

\begin{theorem}\label{thm:main}
For the Milnor fiber $F_f$ given by the germ $f$, we have that
    $$\sigma(F_f)=\sigma(X)+T(f)-C(f),$$
   where $\sigma$ denotes the signature, $X$ is given in \cref{ss:4mfld} and $T,C$ denote the triple point number and cross-cap number of $f$ (see \cref{lemma:TC}).
\end{theorem}
\begin{proof}
    We apply the ESzTS formula \cref{th:saekiszucs} for two slice singular manifolds of $\iota|:\partial F_f\hookrightarrow\Sbb^5$: one for the embedding $$ \iota:F_f\hookrightarrow B^6 $$
   and another for the map given in \cref{ss:constructionh}
   $$h:X\to\shat{B}\cong B^6,$$
   which we know that it is a slice singular manifold of $\iota|$ by \cref{thm:reghom of h}. Hence, $i_b\big(\iota|\big)=i_b\big(h|_{\partial X}\big)$ by \cref{th:saekiszucs}.
   
   On the one hand,
   \begin{equation}
   \begin{aligned}
    i_b\big(\iota|\big)&= \frac{3}{2} \big( \sigma (F_f) - \alpha (\partial F_f)\big) +  \frac{1}{2} \big(3 t(\iota) - 3 \ell(\iota) -
 L\big(\iota|_{\partial F_f}\big)\big) \\
 &= \frac{3}{2} \big( \sigma (F_f) - \alpha (\partial F_f)\big),
   \end{aligned}
   \label{eq:iotaib}
   \end{equation}
   since $\iota$ is an embedding. On the other,
      \begin{equation}
   \begin{aligned}
    i_b\big(h|_{\partial X}\big)&= \frac{3}{2} \big( \sigma (X) - \alpha (\partial X)\big) +  \frac{1}{2} \big(3 t(h) - 3 \ell(h) -
 L\big(h|_{\partial X}\big)\big) \\
 &= \frac{3}{2} \big( \sigma (X) - \alpha (\partial F_f)+ T(f) - C(f) \big),
   \end{aligned}
   \label{eq:hib}
   \end{equation}
   since $\partial X\cong \partial F_f$, $t(h)=T(f)$ by \cref{prop:t=T},  $\ell(h)=C(f)$ by \cref{thm:l=c}, and $L\big(h|_{\partial X}\big)=0$ since $h$ is an embedding on $\partial X$. The result follows from \cref{eq:iotaib,eq:hib}.
\end{proof}

There is a notion of $\eqA$-equivalence in the topological category, we simply assume that the isomorphisms of \cref{eq:defAeq} are homeomorphisms.

\begin{definition}
Two holomorphic map germs $f,f'$ are \textit{topologically equivalent}, or $\eqA_{\Ccal^0}$\textit{-equivalent} if $f=\psi\circ f'\circ\phi$ for some homeomorphisms $\psi,\phi$.
\end{definition}

\begin{corollary}\label{cor:signature top}
    The signature $\sigma(F_f)$ is invariant by $\Ascr_{\Ccal^0}$-equivalence of map germs $f$. 
\end{corollary}
\begin{proof}
It is proven in \cite{Bobadilla2019a} that the numbers $T(f)$ and $C(f)$ are invariant by $\eqA_{\Ccal^0}$-equivalence. 

Furthermore,  it is shown in \cite{Pinter2023} that the vertical indices $\mathfrak{vi}_{f\{i\}}$ (see \cref{def:verticalindex}) are also a topological invariant. Indeed, they give a completely topological definition. Since, by \cref{prop:intersection matrix}, the intersection matrix of $X$ uses the vertical indices and the pair-wise intersection multiplicities of the $D_i$ (which coincide with the linking numbers of $\delta_i$); the signature $\sigma(X)$ is also a topological invariant. The result follows by \cref{thm:main}.
\end{proof} 

\section{Examples}\label{s:ex}

\subsection{Previous comments}

Here we compute the signature of the Milnor fiber for several examples using \cref{th:main}, with a summary in \cref{table:examples}. We study all the simple map germs given in \cite[Theorem 1.1]{Mond1985} (the cross-cap $S_0$ and the families $S_k,B_k,C_k,H_k$), whose values $C(f)$ and $T(f)$ are given in \cite[Table 1]{Mond1987} (recall that one can also use the formulas of \cref{lemma:TC}). We also study other examples with other interesting properties later. The signature $\sigma(X)$ depends on the vertical indices corresponding to the untwisted components, see \cref{prop:intersection matrix}, which can be computed using the description in \cref{sss:vertind}. Furthermore, the boundary of the Milnor fiber is presented for some of these examples in \cite{Nemethi2018} (including all simple germs)  and \cref{th:sumvert} is verified for these examples in \cite{Pinter2023}.

Among the simple germs, only the family $H_k$ has triple points, i.e., $T(H_k)>0$. Hence, the rest (the Whitney umbrella, $S_{k}$, $B_k$ and $C_k$ and $F_4$) are of special type: they are corank--1 germs (i.e., $\rk df_0=1$) with $T(f)=0$. Equivalently, they are in Thom--Boardmann class $\Sigma^{1,0}$. For these germs the computation of the vertical indices is simple, we give it now.

A corank--1 germ $f$ with $T(f)=0$ is called a \textit{fold}, they have the normal form (up to $\Ascr$-equivalence) 
\begin{equation*}
f(u,v) = \big( u, v^2, v p(u,v^2) \big),
\end{equation*}
where $p(u,0)\neq0$. In this case, the equation of the image $\big(\im(f),0\big)$ is 
\begin{equation*}
g(x,y,z)= yp^2(x,y) - z^2 =0,
\end{equation*}
and it is easy to see that the double points of $f$ are given by $\bm{D}=\big\{p(u,v^2)=0\big\}$. See \cite{Mond1985,Marar2009} for details. 

Then, we can choose $\tau(x,y,z)=z$ as a transverse section (cf. \cref{sss:vertind}) and we have \cite{Nemethi2018,Pinter2023}:
    \begin{enumerate}[label={\itshape (\roman*)}]
        \item $\mathfrak{v}^{\tau}_{f\{i\}}=0$ for all $i$,
        \item 
         \[ \lambda^{\tau}_i=-\sum_{k \neq i} D_i \cdot D_k -D_i \cdot \{v=0 \},\]
         \item the vertical indices are
\begin{equation*}
\mathfrak{vi}_{f\{i\}} \coloneqq \left\{ \begin{array}{ccc}
\lambda_i^\tau  + \lambda_{i'}^\tau   & \textnormal{ if } & f\{i\}=\left\{i,i'\right\}, \\
\lambda_i^\tau    & \textnormal{ if } & f\{i\} = \left\{i\right\}. \\
\end{array} \right.
\end{equation*}
\item the cross-cap number is
\begin{equation*}
    C(f)=\bm{D} \cdot \{v=0\}=\frac{\mathcal{O}_{(\C^2,0)}}{\big(v, p(u,v^2)\big)}.
\end{equation*}
    \end{enumerate}

\begin{table}[hbt]
\begin{tabular}{lllcSc}\hline
Name & $f(u,v)$   &  $g(x,y,z)$        & \multicolumn{2}{Sc}{$\sigma(F_f)$}\\ \hline
Cross-cap             & $(u,v^2, uv)$                   & $x^2y-z^2$                         & \multicolumn{2}{Sc}{$-1$}          \\
$S_{k}$, $k \geq 1$   & $( u, v^2, v^3 + u^{k+1} v )$   & $y(y+x^{k+1})^2 - z^2$             & $-k-2$           & $-k-1$         \\
$B_{k} $, $k\geq 2$   & $ ( u, v^2, u^2v + v^{2k+1} ) $ & $ y(x^2+y^k)^2 - z^2$              & $-3$             & $-2$           \\
$C_k $, $k\geq 3$     & $( u, v^2, uv^3 + u^k v )$      & $y(xy+x^k)^2 - z^2$                & $-k-1$           & $-k$           \\
$F_4$                 & $( u, v^2, u^3v + v^5 )$        & $y(x^3+y^2)^2 - z^2$               & \multicolumn{2}{Sc}{$-3$}          \\
$H_k$, $k \geq 2$     & $( u, uv + v^{3k-1}, v^3 )$     & $z^{3k-1} - y^3 + x^3 z + 3 xyz^k$ & \multicolumn{2}{Sc}{$k$}         \\
corank--2             & $(u^2, v^2, u^3 + v^3+ uv)$     & See \cref{ex:corank2}              & \multicolumn{2}{Sc}{$-2$}         \\
Triple point          & See \cref{eq:tr}                & $xyz$                              & \multicolumn{2}{Sc}{$0$}\\ \hline
                      &                                 &                                    & $k=2n+1$         & $k=2n$   
\end{tabular}
\caption{Summary of all the examples in \cref{s:ex}.}
\label{table:examples}
\end{table}

\subsection{Simple germs, and a corank--2 example}

\begin{example}[$S_0$, Whitney umbrella, or cross-cap]\label{ex:s0}
For the stable germ
\begin{equation*}
    f(u,v)=(u,v^2, uv),
\end{equation*}
we have
\begin{equation*}
\begin{matrix}
 p(u,v^2)=u, &\  g(x,y,z)=x^2y-z^2,\\[0.2cm]
    C(f)=1,  & \   T(f)=0.
   \end{matrix}
\end{equation*}

Since $p$ is irreducible in $ \mathcal{O}_{(\C^2, 0)}$, $D=\big\{p(u,v^2)=0 \big\}$ is necessarily a twisted component. Indeed, 
\begin{equation*}
    f(0,v)=(0,v^2, 0)
\end{equation*}
is a double covering of its image (with branching locus the origin). With the transverse section $\tau$ as above, the only vertical index is
\begin{equation*}
    \mathfrak{vi}=\lambda^{\tau}=-D \cdot \{v=0 \} =-1.
\end{equation*}
Nevertheless, this information is not relevant for $\sigma(F_f)$, since it depends only on the vertical indices corresponding to the untwisted components. Indeed, $H_2(X, \Z)$ is trivial (see \cref{eq:HkX}), so $\sigma(X)=0$. Hence,
\begin{equation*}
    \sigma(F_f)=\sigma(X)+T(f)-C(f)=0+0-1=-1.
\end{equation*}

For the Whitney umbrella we can determine the Milnor fiber $F_f$ itself. Here we show directly (without using the argument of our paper) that $F_f$ is a disc bundle over $S^2$ with Euler number $-2$, verifying that $\sigma(F_f)=-1$.
\end{example}

\begin{example}[The family $S_{k}$]\label{ex:sk}
For the family 
$$f(u, v) = ( u, v^2, v^3 + u^{k+1} v ),\  k\geq 1$$
we have 
\begin{equation*}
\begin{matrix}
 p(u,v^2)=v^2+u^{k+1}, &\  g(x, y, z) = y(y+x^{k+1})^2 - z^2,\\[0.2cm]
    C(f)=k+1,  & \   T(f)=0.
   \end{matrix}
\end{equation*}

\textit{\textbf{Odd case.}} If $k=2n-1$, then 
 \[p(u,v^2)=(u^n-iv)(u^n+iv),\]
so $D_1=\{ u^n-iv=0 \}$ and $D_2=\{ u^n+iv=0 \}$ form a couple of untwisted components, mapped to the unique component of $f(\bm{D})$, since a point $(u_0,v_0)\in D_1$ has the same image as $(u_0,-v_0)\in D_2$. We have 
\begin{equation*}
\begin{aligned}
    \lambda_1^{\tau}&=-D_1 \cdot D_2- D_1 \cdot \{v=0 \}=-n-n=-k-1 \textnormal{ and}\\
    \lambda_2^{\tau}&=-D_2 \cdot D_1- D_2 \cdot \{v=0 \}=-n-n=-k-1,
   \end{aligned}
\end{equation*}
therefore 
\begin{equation*}
    \mathfrak{vi}=\mathfrak{vi}_{\{1, 2\}}=\lambda_1^{\tau}+\lambda_2^{\tau}=-2k-2.
\end{equation*}
$H_2(X, \Z)\cong\Z$, generated by the cycle $\mathcal{C}_{\{1, 2\}}$ with self-intersection (recall \cref{prop:intersection matrix})
\begin{equation*}
    \big[\mathcal{C}_{\{1, 2\}}\big] \cdot \big[\mathcal{C}_{\{1, 2\}}\big]=2 D_1 \cdot D_2+\mathfrak{vi}=2n-2k-2=-k-1.
\end{equation*}
Thus, the intersection matrix of $X$ is the $1 \times 1$ matrix $(-k-1)$, its signature is $\sigma(X)=-1$. Hence. the signature of the Milnor fiber is
\begin{equation*}
    \sigma(F_f)=\sigma(X)+T(f)-C(f)=-1+0-(k+1)=-k-2.
\end{equation*}

\textit{\textbf{Even case.}} If $k=2n$, then $D$ is irreducible. In that case, there is only one twisted component, which is irrelevant for the signature. $H_2(X, \Z)$ is trivial, so $\sigma(X)=0$. Hence,
\begin{equation*}
    \sigma(F_f)=\sigma(X)+T(f)-C(f)=0+0-(k+1)=-k-1.
\end{equation*}  
\end{example}

\begin{example}[The family $B_{k}$]
For the family 
$$f(u, v) = ( u, v^2, u^2v + v^{2k+1} ),\  k\geq 2$$
we have 
\begin{equation*}
\begin{matrix}
 p(u,v^2)=u^2+v^{2k}, &\  g(x, y, z) = y(x^2+y^k)^2 - z^2,\\[0.2cm]
    C(f)=2,  & \   T(f)=0.
   \end{matrix}
\end{equation*}
 
 Although $\bm{D}$ always has two components $D_1=\{u+iv^k=0\}$ and $D_2=\{u-iv^k=0\}$, whether they are twisted or untwisted, or the number of the components of $f(\bm{D})$, depends on the parity of $k$.
\newline

 \textbf{Odd case.} If $k=2n+1$, any point $(u_0,v_0)\in D_1$ has the same image as $(u_0,-v_0)\in D_2$. Therefore, $D_1$ and $D_2$ are a pair of components corresponding to the unique untwisted component of $f(\bm{D})$.
 Then, 
$$\begin{aligned}
\lambda_1^{\tau}=\lambda_2^{\tau}=&-D_1 \cdot D_2-D_i \cdot \{v=0\}=-k-1,\text{ and}\\
 \mathfrak{vi}_{\{1, 2\}}=&\lambda_1^{\tau}+\lambda_2^{\tau}=-2k-2.\end{aligned}$$
$H_2(X)\cong\Z$ and it is generated by the cycle $\mathcal{C}_{\{1, 2\}}$. Its self-intersection is 
\begin{equation*}
    \big[\mathcal{C}_{\{1, 2\}}\big] \cdot \big[\mathcal{C}_{\{1, 2\}}\big]=2 D_1 \cdot D_2+\mathfrak{vi}=2k+(-2k-2)=-2.
\end{equation*}
Hence $\sigma(X)=-1$. The signature of the Milnor fiber is
\begin{equation}
    \sigma(F_f)=\sigma(X)+T(f)-C(f)=-1+0-2=-3.
\end{equation}


\textit{\textbf{Even case.}} If $k=2n$, then $(u_1,v_1),(u_1,-v_1)\in D_1$ have the same image and $(u_2,v_2),(u_2,-v_2)\in D_2$ as well, so $f(\bm{D})$ has two twisted components corresponding to $D_1$ and $D_2$. $H_2(X)$ is trivial, so $\sigma(X)=0$ and the signature of the Milnor fiber is
\begin{equation*}
    \sigma(F_f)=\sigma(X)+T(f)-C(f)=0+0-2=-2.
\end{equation*}
%

\end{example}

\begin{example}[The family $C_k $]
For the family 
$$f(u, v) = ( u, v^2, uv^3 + u^k v ),\  k\geq 3 $$
we have 
\begin{equation*}
\begin{matrix}
 p(u,v^2)=uv^2+u^k, &\  g(x, y, z) = y(xy+x^k)^2 - z^2,\\[0.2cm]
    C(f)=k,  & \   T(f)=0.
   \end{matrix}
\end{equation*}

\textbf{\textit{Odd case.}} If $k=2n+1$, then
$$  p(u,v^2)=uv^2+u^k= u(u^n+iv)(u^n-iv).$$
It is easy to see that $(0,v_1),(0,-v_1)\in D_1=\left\{u=0\right\}$ are mapped to the same point, so $f(D_1)$ is a twisted component. Similarly, $(u_0,v_0)\in D_2=\left\{u^n+iv=0\right\}$ and 
$(u_0,-v_0)\in D_3=\left\{u^n-iv=0\right\}$ are mapped to the same point, so $f(D_2)=f(D_3)$ is an untwisted component. Since only untwisted components are relevant, we compute
\[\begin{aligned} 
\lambda_2^{\tau}=\lambda_3^{\tau}=&-D_2 \cdot D_1-D_2 \cdot D_3-D_2 \cdot \{v=0\}=-1-n-n=-k\text{ and}\\
\mathfrak{vi}_{\left\{2,3\right\}}=&-2k
\end{aligned} \]
$H_2(X)\cong\Z$ is generated by $\mathcal{C}_{\{2, 3\}}$ with self-intersection
\begin{equation*}
    \big[\mathcal{C}_{\{2, 3\}}\big] \cdot \big[\mathcal{C}_{\{2, 3\}}\big]=2 D_2 \cdot D_3+\mathfrak{vi}_B=2n-2k=-k-1.
\end{equation*}
Therefore, $\sigma(X)=-1$ and the signature of the Milnor fiber is
\begin{equation*}
    \sigma(F_f)=\sigma(X)+T(f)-C(f)=-1+0-k=-k-1.
\end{equation*}

%
 \textbf{Even case.} If $k=2n$,then
$$  p(u,v^2)=uv^2+u^k= u(v^2+u^{2n-1}).$$
It is easy to see that $D_1=\left\{u=0\right\}$ and $D_2=\left\{v^2+u^{2n-1}=0\right\}$ cannot have common images away from the origin, so they must give two twisted components $f(D_1)$ and $f(D_2)$. Since $H_2(X)$ is trivial, the signature of the Milnor fiber is
\begin{equation*}
    \sigma(F_f)=\sigma(X)+T(f)-C(f)=0+0-k=-k.
\end{equation*}


\end{example}

\begin{example}[The germ $F_4$] For the germ
     \begin{equation*}
         f(u, v) = ( u, v^2, u^3v + v^5 ) 
     \end{equation*} 
   we have 
\begin{equation*}
\begin{matrix}
 p(u,v^2)=u^3+v^4, &\  g(x, y, z) = y(x^3+y^2)^2 - z^2,\\[0.2cm]
    C(f)=3,  & \   T(f)=0.
   \end{matrix}
\end{equation*}
%

Since $D$ is irreducible, it defines one twisted component and $H_2(X)$ is trivial, so $\sigma(X)=0$. The signature of the Milnor fiber is 
\begin{equation*}
    \sigma(F_f)= \sigma(X)+T(f)-C(f)=-3.
\end{equation*}
\end{example}

The following example is the unique family of simple germs that is not a fold singularity, so $T(f)>0$. The usual simplifications for fold singularities are not valid any more, however, the vertical indices can be computed directly according to \cref{sss:vertind}.

\begin{example}[The family $H_{k} $]
Let us consider the family
$$f(u, v) = ( u, uv + v^{3k-1}, v^3 ),\  k\geq 2. $$

In this case, the equation of the image $\big(\im(f),0\big)$ is
$$ g(x, y, z) = z^{3k-1} - y^3 + x^3 z + 3 xyz^k; $$
the equation of $\bm{D}$ is
$$ (u-\rho v^{3k-2})(u-\rho^2 v^{3k-2}),$$
where $\rho=e^{\frac{i2\pi}{3}}$; and
$$    C(f)=2,   \   T(f)=k-1.$$


The equation $g$ and $T(f)$ can be computed by using fitting ideals (see, e.g., \cite[Example 1.3.4.]{Gergothesis}). We observe that $(u_0,v_0)\in D_1=\big\{u-\rho v^{3k-2}=0\big\}$ and $(u_0,\rho v_0)\in D_2=\big\{u-\rho^2 v^{3k-2}=0\big\}$ have the same image. Hence, $f(\bm{D})$ has only one untwisted component.

The computation of the vertical index is much more complicated than in the previous examples but it follows the general method (cf. \cref{sss:vertind}), it was done in \cite[Section 6.7]{Nemethi2018}, and the result is $\mathfrak{vi}=-3k-1$. 
$H_2(X)\cong\Z$  is generated by   $ \mathcal{C}_{\{1, 2\}}$ with self-intersection 
\begin{equation*}
    \big[\mathcal{C}_{\{1, 2\}}\big] \cdot \big[\mathcal{C}_{\{1, 2\}}\big]=2 D_1 \cdot D_2+\mathfrak{vi}=2\cdot (3k-2)-3k-1=3k-5.
\end{equation*}
The intersection matrix of $X$ is $(3k-5)$. Hence 
\begin{equation*}
    \sigma(X)=1
\end{equation*}
since $k\geq2$.
The signature of the Milnor fiber is
\begin{equation*}
    \sigma(F_f)=\sigma(X)+T(f)-C(f)=k .
\end{equation*}
\end{example}

\begin{note}
The previous example provides a family of non-isolated singularities such that their Milnor fiber have signature arbitrarily big. This could be interesting, for example, when one compares this example with Durfee's conjecture for isolated singularities \cite{Durfee1978} (cf. \cite{Wahl1981,Kollar2017}).
\end{note}

The previous examples have corank one, i.e., $\rk df_0=1$. This condition usually simplifies other kind of computations (e.g., \cite{Goryunov1993,Houston2010,GimenezConejero2022b}). We present now the computation on a corank two germ, $\rk df_0=0$, from \cite{Marar2008}.

\begin{example}[A corank two germ]\label{ex:corank2}
For the germ
\[ f (s, t) = (u^2, v^2, u^3 + v^3+ uv) ,\]
the equation of the image is
\[g(x, y, z) = x^6-2x^4y + x^2y^2 - 2x^3y^3-2xy^4 + y^6 - 8x^2y^2z - 2x^3z^2-2xyz^2-2y^3z^2+z^4  ,\]
and the equation of the double points $\bm{D}$ is
\[(u+v^2)(u^2+v)(u+v)(u+\rho v)(u+\rho^{2}v)\]
where $\rho=e^{\frac{i2\pi}{3}}$ (see \cite{Marar2008}).
It is also easy to compute with \cref{lemma:TC} that
\[C(f)=3, \ T(f)=1.\]

All the five components of $f(\bm{D})$ are twisted, since the following pairs of points have the same image
$$\begin{aligned}
(u_1,v_1),(u_1,-v_1)\in D_1&=\left\{u+v^2=0\right\},\\
(u_1,v_1),(-u_1,v_1)\in D_2&=\left\{u^2+v=0\right\},\\
(u_1,v_1),(-u_1,-v_1)\in D_3&=\left\{u+v=0\right\},\\
(u_1,v_1),(-u_1,-v_1)\in D_4&=\left\{u+\rho v=0\right\},\\
(u_1,v_1),(-u_1,-v_1)\in D_5&=\left\{u+\rho^{2}v=0\right\}.\\
\end{aligned}$$
Since all are twisted, $H_2(X)$ is trivial and $\sigma(X)=0$. 

The signature of the Milnor fiber is
\begin{equation}
    \sigma(F_f)=\sigma(X)+T(f)-C(f)=0+1-3=-2.
\end{equation}

 \end{example}

\subsection{A multigerm example}

%
%

   One might guess that the terms $-C(f)$ and $+T(f)$ in the signature formula \cref{th:main} correspond to the signature of the Milnor fiber of the cross-cap points and the triple values  respectively, which is really $-1$ for a cross-cap, and it is expected to be $+1$ for a triple value. Indeed, it is true that the Milnor fiber $F_f$ of $\big(\im(f),0\big)$ coincides with the nearby fiber of the image of a stable deformation $f_s$ (by an argument of stratification theory using the stratified version of Ehresmann's lemma, \cite[Theorem 5.2]{Gibson1976a}), and the latter one might be patched together from the Milnor fibers of the cross caps and the triple values (this idea was first used in a preliminary version of \cite{Bobadilla2019a}). One would like to use additivity of the signature to finish the argument. 
   
   However, the signature of the Milnor fiber of a triple value is 0, not 1. This is reflecting the fact that one also has to take into account the double points that connect the different cross-caps and triple points, changing the signature. This is, in some sense, also what the term $\sigma(X)$ does in our formula.

   We proceed to, first, determine the Milnor fiber of a triple value directly, and show that its signature is 0. Then, we verify it by using our argument. Since the triple value is a multigerm, \cref{th:main} cannot be applied in the presented form, but we can use the same method as we deduced the signature formula for monogerms to compute the signature of the Milnor fiber for a triple value. This example can be considered as the first step toward the generalization of \cref{th:main} for multigerms.

 A triple value, more precisely, a normal crossing intersection of three regular branches, is the $\mathscr{A}$-equivalence class of the multigerm
    \begin{equation}\label{eq:tr}
    f: \bigsqcup_{i=1}^3 (\C^2,0)_i \to (\C^3, 0), \
    \left\{ \begin{array}{ccc} 
     f(u_1,v_1)    & = & (0,u_1,v_1)  \\
     f(u_2,v_2)    & = & (u_2,0,v_2)  \\
     f(u_3,v_3)    & = & (u_3,v_3,0)  \\
    \end{array} \right.
\end{equation}
with image $g(x,y,z)=xyz=0$. The radius of the Milnor ball can be chosen to be $\epsilon =1$, and the Milnor fiber $F_f$ of the triple value is defined by the system
\begin{equation*}
    \left.
\begin{array}{ccc}
|x|^2+|y|^2+|z|^2 & \leq & 1 \\
xyz & = & \delta
\end{array}
    \right\}
\end{equation*}
where $0 < \delta \ll 1$. By projecting to the $x$ coordinate one can see that $F_f$ is a locally trivial bundle of cylinders $\Sbb^1 \times [-1,1]$ over the base space $\Sbb^1 \times [-1,1]$. Since $F_f$ is oriented, it could be either the trivial bundle or $Z_{\partial} \times [-1,1]$ (recall \cref{eq:Y}), but the middle torus 
\begin{equation*}
    \Tbb^2=\{(x,y,z) \in F_f \ | \ |x|=|y|=|z|=\sqrt[3]{\delta}\} 
\end{equation*}
shows that $F_f$ is indeed the trivial bundle,
\begin{equation*}
    F_f \cong \Sbb^1 \times \Sbb^1 \times [-1,1] \times [-1,1] \cong \Tbb^2 \times B^2.
\end{equation*}

Hence, $H_2(F_f, \Z)\cong\Z$ is generated by the cycle $\Tbb^2 \times \{0\}$. Since its self-intersection is 0, the intersection matrix of $F_f$ is $(0)$ and $\sigma(F_f)=0$.

\begin{example}[ESzTS formula for the Milnor fiber of a triple value]
Although both our signature formula \cref{th:main} and the boundary formula \cref{th:milnperpart1} is formulated for monogerms, we can follow the same method  to verify  $\sigma(F_f)=0$ for a triple value. Using the double points of $f$ in the source,  we construct a $4$-manifold $X$ with boundary $\partial X \cong \partial F_f$, and a slice singular manifold $h: X \to B^6$ for the inclusion $\partial F \hookrightarrow \Sbb^5$, then we compute $\sigma(F_f)$ by \cref{th:saekiszucs}.

The set of double values is
  $$\{(x,0,0)\} \cup \{(0,y,0)\} \cup \{(0,0,z)\} \subset B^6 \subset \C^3,$$
with preimages in each copy of $(\CC^2,0)$:
  $$D_i\coloneqq  D_{u,i} \cup D_{v,i}\coloneqq\{v_i=0\} \cup \{u_i=0\} \subset (\C^2,0)_i$$
for $i\in\{1,2,3\}$. Note that the link $\delta_i=\delta_{u,i}\sqcup\delta_{v,i}$ of $D_i$ is a Hopf link. We also have the pairings
\begin{equation*}
    f(D_{u,1}) = f(D_{v,3}) ,\ 
    f(D_{u,2}) = f(D_{u,3}), \ 
    f(D_{v,1})=f(D_{v,2}).
\end{equation*}

Let $N_{u,i}$ and $N_{v,i}$ be tubular neighbourhoods of $\delta_{u,i}$ and $\delta_{v,i}$ respectively. 
Following \cref{eq:def of X,eq:def of X2}, we define
\begin{equation*}
    X=(B^4_1 \sqcup B_2^4 \sqcup B^4_3  ) / \Phi
\end{equation*}
with the gluing of the related tori
\begin{equation*}
    \Phi=\left\{
\begin{array}{ccc}
     N_{u,1} & \to & - N_{v,3} \\
       N_{u,2} & \to & -N_{u,3} \\
         N_{v,1} & \to & - N_{v,2} \\
\end{array}
    \right.
\end{equation*}
in the trivial way, i.e., by the torus
\begin{equation*}
	\left( \begin{array}{cc}
	-1 & 0 \\
	0 & 1 \\
	\end{array} \right),
	\end{equation*} 
 using the trivialization $(\textnormal{meridian}, \textnormal{topological longitude})$ of the boundaries. In other words, the vertical indices corresponding to different branches are 0.

 Obviously, $X$ can be smoothen. Its boundary is $\partial X \cong \partial F_f \cong \Sbb^1 \times \Sbb^1 \times \Sbb^1$. The second homology of $X$ is $H_2(X, \Z)\cong \Z \oplus \Z \oplus \Z$ generated by the cycles 
 \begin{equation*}
\mathcal{C}_y=D_{u,1} \cup_{\Phi} D_{v,3}, \
\mathcal{C}_x=D_{u,2} \cup_{\Phi} D_{u,3}, \
\mathcal{C}_z=D_{v,1} \cup_{\Phi} D_{v,2}.
 \end{equation*}
 Each generator has zero self-intersection, because of the trivial gluing, and each pairwise intersection is $1$. Hence, the intersection matrix of $X$ is
 \begin{equation*}
	\left( \begin{array}{ccc}
	0 & 1 & 1 \\
	1 & 0 & 1 \\
 1 & 1 & 0 \\
	\end{array} \right),
	\end{equation*} 
 whose eigenvalues are $2, -1, -1$. Therefore, $\sigma(X)=-1$.
 
 For the definition of the map $h: X \to B^6 $, the gluing $\Phi$ is replaced by the gluing of three pieces of $\Sbb^1 \times B^2 \times [-1,1]$, the same way as we did in \cref{rem:othergluing}. Since in this case $f$ is stable both in the holomorphic and in the real $\mathcal{C}^{\infty}$ sense, $f^{\R}_s=f_s=f$. Its extension $h: X \to \shat{B} \cong B^6$ is constructed in the same way we did it in \cref{ss:constructionh}, by creating three one-parameter families of generalized real Whitney umbrellas along $\Sbb^1 \times \{0\} \times \{0\} \subset \Sbb^1 \times B^2 \times [-1,1]$. The linking numbers $\ell$ and $L$ are zero, similarly to \cref{ssec:properties}.

 Then, applying the formula of \cref{th:saekiszucs} to $F_f \hookrightarrow B^6$ and to $h$, 
 \begin{equation*}
     \sigma(F_f)=\sigma(X)-\ell(h)+t(f)=-1+0+1=0.
 \end{equation*}

\end{example}

\bibliographystyle{C:/Users/rgc19/AppData/Local/Programs/MiKTeX/bibtex/bst/bibtex/myalpha.bst}
\bibliography{C:/Users/rgc19/Documents/Textos-Matematicas/Bibliografias/bibtex/bib/mybibs/FullBib.bib} 

\begin{thebibliography}{MNnBPnS12}

\bibitem[dBMN14]{Bobadilla2014}
J.~F. de~Bobadilla and A.~Menegon~Neto.
\newblock The boundary of the {M}ilnor fibre of complex and real analytic
  non-isolated singularities.
\newblock {\em Geometriae Dedicata}, 173:143--162, 2014.

\bibitem[Dur78]{Durfee1978}
A.~H. Durfee.
\newblock The signature of smoothings of complex surface singularities.
\newblock {\em Mathematische Annalen}, 232(1):85--98, 1978.

\bibitem[Ekh01]{Ekholm2001}
T.~Ekholm.
\newblock Differential 3-knots in 5-space with and without self-intersections.
\newblock {\em Topology. An International Journal of Mathematics},
  40(1):157--196, 2001.

\bibitem[ES03]{Ekholm2003}
T.~Ekholm and A.~Sz\"{u}cs.
\newblock Geometric formulas for {S}male invariants of codimension two
  immersions.
\newblock {\em Topology. An International Journal of Mathematics},
  42(1):171--196, 2003.

\bibitem[ES06]{Ekholm2006}
T.~Ekholm and A.~Sz\H{u}cs.
\newblock The group of immersions of homotopy {$(4k-1)$}-spheres.
\newblock {\em The Bulletin of the London Mathematical Society},
  38(1):163--176, 2006.

\bibitem[ET11]{Ekholm2011}
T.~Ekholm and M.~Takase.
\newblock Singular {S}eifert surfaces and {S}male invariants for a family of
  3-sphere immersions.
\newblock {\em Bulletin of the London Mathematical Society}, 43(2):251--266,
  2011.

\bibitem[FdBPnSS22]{Bobadilla2019a}
J.~Fern\'{a}ndez~de Bobadilla, G.~Pe\~{n}afort Sanchis, and E.~Sampaio.
\newblock Topological invariants and {M}ilnor fibre for $\mathcal{A}$-finite
  germs $\mathbb{C}^2\to\mathbb{C}^3$.
\newblock {\em Dalat University Journal of Science}, Volume 12(Issue 2), 2022.

\bibitem[FQ90]{Freedman1990}
M.~H. Freedman and F.~Quinn.
\newblock {\em Topology of 4-manifolds}, volume~39 of {\em Princeton
  Mathematical Series}.
\newblock Princeton University Press, Princeton, NJ, 1990.

\bibitem[Fre82]{Freedman1982}
M.~H. Freedman.
\newblock The topology of four-dimensional manifolds.
\newblock {\em Journal of Differential Geometry}, 17(3):357--453, 1982.

\bibitem[GC21]{Robertothesis}
R.~Gim{\'e}nez~Conejero.
\newblock {\em Singularities of germs and vanishing homology}.
\newblock PhD thesis, Universitat de Val{\`e}ncia, 2021.

\bibitem[GCNB22]{GimenezConejero2022b}
R.~Gim\'enez~Conejero and J.~J. Nu{\~{n}}o-Ballesteros.
\newblock Singularities of mappings on icis and applications to whitney
  equisingularity.
\newblock {\em Advances in Mathematics}, 408:108660, 2022.

\bibitem[GM93]{Goryunov1993}
V.~Goryunov and D.~Mond.
\newblock Vanishing cohomology of singularities of mappings.
\newblock {\em Compositio Mathematica}, 89(1):45--80, 1993.

\bibitem[GS99]{Gompf1999}
R.~E. Gompf and A.~I. Stipsicz.
\newblock {\em {$4$}-manifolds and {K}irby calculus}, volume~20 of {\em
  Graduate Studies in Mathematics}.
\newblock American Mathematical Society, Providence, RI, 1999.

\bibitem[GWdPL76]{Gibson1976a}
C.~G. Gibson, K.~Wirthm\"{u}ller, A.~A. du~Plessis, and E.~J.~N. Looijenga.
\newblock {\em Topological stability of smooth mappings}.
\newblock Lecture Notes in Mathematics, Vol. 552. Springer-Verlag, Berlin-New
  York, 1976.

\bibitem[HM85]{Hughes1985}
J.~F. Hughes and P.~M. Melvin.
\newblock The {S}male invariant of a knot.
\newblock {\em Commentarii Mathematici Helvetici}, 60(4):615--627, 1985.

\bibitem[Hou10]{Houston2010}
K.~Houston.
\newblock Stratification of unfoldings of corank 1 singularities.
\newblock {\em The Quarterly Journal of Mathematics}, 61(4):413--435, 2010.

\bibitem[Juh05]{Juhasz2005}
A.~Juh\'{a}sz.
\newblock A geometric classification of immersions of 3-manifolds into 5-space.
\newblock {\em Manuscripta Mathematica}, 117(1):65--83, 2005.

\bibitem[Kin78]{King1978}
H.~C. King.
\newblock Topological type of isolated critical points.
\newblock {\em Annals of Mathematics. Second Series}, 107(2):385--397, 1978.

\bibitem[Kin15]{Kinjo2015}
S.~Kinjo.
\newblock Immersions of 3-sphere into 4-space associated with {D}ynkin diagrams
  of types {$A$} and {$D$}.
\newblock {\em Bulletin of the London Mathematical Society}, 47(4):651--662,
  2015.

\bibitem[KN17]{Kollar2017}
J.~Koll\'{a}r and A.~N\'{e}methi.
\newblock Durfee's conjecture on the signature of smoothings of surface
  singularities.
\newblock {\em Annales Scientifiques de l'\'{E}cole Normale Sup\'{e}rieure.
  Quatri\`eme S\'{e}rie}, 50(3):787--798, 2017.
\newblock With an appendix by Tommaso de Fernex.

\bibitem[KNS14]{Katanaga2014}
A.~Katanaga, A.~N\'{e}methi, and A.~Sz\H{u}cs.
\newblock Links of singularities up to regular homotopy.
\newblock {\em Journal of Singularities}, 10:174--182, 2014.

\bibitem[Lau77]{Laufer1977}
H.~B. Laufer.
\newblock On {$\mu $} for surface singularities.
\newblock In {\em Several complex variables ({P}roc. {S}ympos. {P}ure {M}ath.,
  {V}ol. {XXX}, {P}art 1, {W}illiams {C}oll., {W}illiamstown, {M}ass., 1975)},
  pages 45--49. Amer. Math. Soc., Providence, R.I., 1977.

\bibitem[L{\^e}73]{Le1973}
D.~T. L{\^e}.
\newblock Calcul du nombre de cycles \'{e}vanouissants d'une hypersurface
  complexe.
\newblock {\em Ann. Inst. Fourier (Grenoble)}, 23(4):261--270, 1973.

\bibitem[Mil68]{Milnor1968}
J.~Milnor.
\newblock {\em Singular points of complex hypersurfaces}.
\newblock Annals of Mathematics Studies, No. 61. Princeton University Press,
  Princeton, N.J.; University of Tokyo Press, Tokyo, 1968.

\bibitem[MM89]{Marar1989}
W.~L. Marar and D.~Mond.
\newblock Multiple point schemes for corank {$1$} maps.
\newblock {\em Journal of the London Mathematical Society. Second Series},
  39(3):553--567, 1989.

\bibitem[MNB20]{Mond2020}
D.~Mond and J.~J. Nu{\~{n}}o-Ballesteros.
\newblock {\em Singularities of mappings}, volume 357 of {\em Grundlehren der
  mathematischen Wissenschaften}.
\newblock Springer, Cham, 2020.

\bibitem[MNnB08]{Marar2008}
W.~L. Marar and J.~J. Nu\~{n}o Ballesteros.
\newblock A note on finite determinacy for corank 2 map germs from surfaces to
  3-space.
\newblock {\em Mathematical Proceedings of the Cambridge Philosophical
  Society}, 145(1):153--163, 2008.

\bibitem[MNnB09]{Marar2009}
W.~L. Marar and J.~J. Nu\~{n}o Ballesteros.
\newblock The doodle of a finitely determined map germ from {$\Bbb R^2$} to
  {$\Bbb R^3$}.
\newblock {\em Advances in Mathematics}, 221(4):1281--1301, 2009.

\bibitem[MNnB14]{Marar2014}
W.~L. Marar and J.~J. Nu\~{n}o Ballesteros.
\newblock Slicing corank 1 map germs from {$\mathbb{C}^2$} to {$\mathbb{C}^3$}.
\newblock {\em The Quarterly Journal of Mathematics}, 65(4):1375--1395, 2014.

\bibitem[MNnBPnS12]{Marar2012}
W.~L. Marar, J.~J. Nu\~{n}o Ballesteros, and G.~Pe\~{n}afort Sanchis.
\newblock Double point curves for corank 2 map germs from {$\mathbb{C}^2$} to
  {$\mathbb{C}^3$}.
\newblock {\em Topology and its Applications}, 159(2):526--536, 2012.

\bibitem[Mon85]{Mond1985}
D.~Mond.
\newblock On the classification of germs of maps from {${\bf R}^2$} to {${\bf
  R}^3$}.
\newblock {\em Proceedings of the London Mathematical Society. Third Series},
  50(2):333--369, 1985.

\bibitem[Mon87]{Mond1987}
D.~Mond.
\newblock Some remarks on the geometry and classification of germs of maps from
  surfaces to {$3$}-space.
\newblock {\em Topology}, 26(3):361--383, 1987.

\bibitem[Mon91]{Mond1991}
D.~Mond.
\newblock Vanishing cycles for analytic maps.
\newblock In {\em Singularity theory and its applications, {P}art {I}
  ({C}oventry, 1988/1989)}, volume 1462 of {\em Lecture Notes in Math.}, pages
  221--234. Springer, Berlin, 1991.

\bibitem[MPW07]{Michel2007}
F.~Michel, A.~Pichon, and C.~Weber.
\newblock The boundary of the {M}ilnor fiber of {H}irzebruch surface
  singularities.
\newblock In {\em Singularity theory}, pages 745--760. World Sci. Publ.,
  Hackensack, NJ, 2007.

\bibitem[MS92]{Massey1992}
D.~B. Massey and D.~Siersma.
\newblock Deformation of polar methods.
\newblock {\em Universit\'{e} de Grenoble. Annales de l'Institut Fourier},
  42(4):737--778, 1992.

\bibitem[N\'99]{Nemethi1999}
A.~N\'{e}methi.
\newblock Five lectures on normal surface singularities.
\newblock In {\em Low dimensional topology ({E}ger, 1996/{B}udapest, 1998)},
  volume~8 of {\em Bolyai Soc. Math. Stud.}, pages 269--351. J\'{a}nos Bolyai
  Math. Soc., Budapest, 1999.
\newblock With the assistance of \'{A}gnes Szil\'{a}rd and S\'{a}ndor
  Kov\'{a}cs.

\bibitem[NP15]{Nemethi2015}
A.~N\'{e}methi and G.~Pint\'{e}r.
\newblock Immersions associated with holomorphic germs.
\newblock {\em Commentarii Mathematici Helvetici. A Journal of the Swiss
  Mathematical Society}, 90(3):513--541, 2015.

\bibitem[NP18]{Nemethi2018}
A.~N\'{e}methi and G.~Pint\'{e}r.
\newblock The boundary of the {M}ilnor fibre of certain non-isolated
  singularities.
\newblock {\em Periodica Mathematica Hungarica. Journal of the J\'{a}nos Bolyai
  Mathematical Society}, 77(1):34--57, 2018.

\bibitem[NS12]{Nemethi2012}
A.~N\'{e}methi and A.~Szil\'{a}rd.
\newblock {\em Milnor fiber boundary of a non-isolated surface singularity},
  volume 2037 of {\em Lecture Notes in Mathematics}.
\newblock Springer, Heidelberg, 2012.

\bibitem[OS03]{Ozsvath2003}
P.~Ozsv\'{a}th and Z.~Szab\'{o}.
\newblock Knot {F}loer homology and the four-ball genus.
\newblock {\em Geometry and Topology}, 7:615--639, 2003.

\bibitem[Per85]{Perron1985}
B.~Perron.
\newblock Conjugaison topologique des germes de fonctions holomorphes \`a
  singularit\'{e} isol\'{e}e en dimension trois.
\newblock {\em Inventiones Mathematicae}, 82(1):27--35, 1985.

\bibitem[Pin18]{Gergothesis}
G.~Pint{\'{e}}r.
\newblock On certain complex surface singularities, ph.d. thesis.
\newblock {\em Eötvös Loránd University, arXiv:1904.12778}, 2018.

\bibitem[PS23]{Pinter2023a}
G.~Pintér and A.~Sándor.
\newblock Cross-caps, triple points and a linking invariant for finitely
  determined germs.
\newblock {\em Revista Matemática Complutense}, February 2023.

\bibitem[PT23]{Pinter2023}
G.~Pintér and T.~Terpai.
\newblock The boundary of the milnor fibre and a linking invariant of finitely
  determined germs.
\newblock {\em arXiv:2304.12672}, April 2023.

\bibitem[Sae89]{Saeki1989}
O.~Saeki.
\newblock Topological types of complex isolated hypersurface singularities.
\newblock {\em Kodai Mathematical Journal}, 12(1):23--29, 1989.

\bibitem[Sie91a]{Siersma1991}
D.~Siersma.
\newblock Vanishing cycles and special fibres.
\newblock In {\em Singularity theory and its applications, {P}art {I}
  ({C}oventry, 1988/1989)}, volume 1462 of {\em Lecture Notes in Math.}, pages
  292--301. Springer, Berlin, 1991.

\bibitem[Sie91b]{Siersma1991a}
D.~Siersma.
\newblock Variation mappings on singularities with a {$1$}-dimensional critical
  locus.
\newblock {\em Topology. An International Journal of Mathematics},
  30(3):445--469, 1991.

\bibitem[Sma59]{Smale1959}
S.~Smale.
\newblock The classification of immersions of spheres in {E}uclidean spaces.
\newblock {\em Annals of Mathematics. Second Series}, 69:327--344, 1959.

\bibitem[SST02]{Saeki2002}
O.~Saeki, A.~Sz{\H{u}}cs, and M.~Takase.
\newblock Regular homotopy classes of immersions of 3-manifolds into 5-space.
\newblock {\em Manuscripta Mathematica}, 108(1):13--32, 2002.

\bibitem[Tak07]{Takase2007}
M.~Takase.
\newblock An {E}kholm-{S}z{\H{u}}cs-type formula for codimension one immersions
  of 3-manifolds up to bordism.
\newblock {\em Bulletin of the London Mathematical Society}, 39(1):39--45,
  2007.

\bibitem[Wah81]{Wahl1981}
J.~Wahl.
\newblock Smoothings of normal surface singularities.
\newblock {\em Topology. An International Journal of Mathematics},
  20(3):219--246, 1981.

\end{thebibliography}
\end{document}